\documentclass[aos,MSNbibl,seceqn,nameyear,dvips]{arximspdf}
\usepackage{dcolumn}
\usepackage{graphicx}

\doi{10.1214/12-AOS977} %kopijuoti is PTS
\volume{40}
\issue{2}
\pubyear{2012}
\firstpage{759}
\lastpage{784}

\makeatletter
\newcolumntype{d}[1]{D{.}{.}{#1}}
\newtheorem{theorem}{Theorem}
\newtheorem{corollary}{Corollary}
\newtheorem{lemma}{Lemma}
\newtheorem{proposition}{Proposition}
\newproclaim{assumption}{Assumption}
\newproclaim{remark}{Remark}
\makeatother

\begin{document}
\begin{frontmatter}

\title{Modeling high-frequency financial data by pure jump processes}
\runtitle{Pure jump modeling}

\begin{aug}
\author[A]{\fnms{Bing-Yi} \snm{Jing}\ead
[label=e1]{majing@ust.hk}\thanksref{t1}},
\author[B]{\fnms{Xin-Bing} \snm{Kong}\corref{}\ead[label=e2]{kongxb@fudan.edu.cn}\thanksref{t2}}
\and
\author[C]{\fnms{Zhi} \snm{Liu}\ead[label=e3]{liuzhi@ust.hk}}
\thankstext{t1}{Supported in part by HK RGC Grants HKUST6011/07P,
HKUST6015/08P, and HKUST6019/10P.}
\thankstext{t2}{Supported in part by the Humanity and Social Science Youth Foundation
of Chinese Ministry of Education No. 12YJC910003.}
\runauthor{B.-Y. Jing, X.-B. Kong and Z. Liu}
\affiliation{Hong Kong University of Science and Technology, Fudan University
and Xiamen University}
\address[A]{B.-Y. Jing\\
Department of Mathematics\\
Hong Kong University of Science and Technology\\
Clear Water Bay, Kowloon\\
Hong Kong\\
\printead{e1}} %adresu isvedimo komanda gale!
\address[B]{X.-B. Kong\\
R636\\
Siyuan Building\\
Guoshun Road 670, Shanghai\\
P. R. China\\
\printead{e2}}
\address[C]{Z. Liu\\
Flat B406\\
Economic Building\\
Xiamen University\\
Xiamen, Fujian Province\\
P. R. China\\
\printead{e3}}
\end{aug}

% HISTORY:
\received{\smonth{1} \syear{2011}}
\revised{\smonth{9} \syear{2011}}

% ABSTRACT
%
\begin{abstract}
It is generally accepted that the asset price processes contain
jumps. In fact, pure jump models have been widely used to model
asset prices and/or stochastic volatilities.
% by pure jump processes with no diffusion components.
The question is: is there any statistical evidence from the
high-frequency financial data to support using pure jump models alone?
The purpose of this paper is to develop such a statistical test
against the necessity of a diffusion component. The test is very
simple to use and yet effective. Asymptotic properties of the
proposed test statistic will be studied.
% The advantages of the test over some alternative ones in the
%literature
% will be discussed.
Simulation studies and some real-life examples are included to
illustrate our results.
\end{abstract}

% KEYWORDS
%

\end{frontmatter}

%s1 ###
\section{Introduction}\label{sec1}

It is now widely accepted that the asset price processes contain
jumps. This is partially based on many empirical evidences, such as
heavy tails in the asset returns;
% which strongly suggest the existence of jumps;
see Cont and Tankov (\citeyear{ConTan04}) and Carr et al. (\citeyear{Caretal02}) and
references therein. In the meantime, many statistical tests have
been established to detect jumps from discretely observed prices [e.g.,
%% To name a few, they are
Jiang and Oomen (\citeyear{JiaOom}), Barndorff-Neilsen and Shepard (\citeyear{BarShe06}), Lee
and Mykland (\citeyear{LeeMyk07}), A\"{\i}t-Sahalia and Jacod (\citeyear{AtSJac10})], and these
test results all seem to support the claim of the existence of jumps
for the asset returns under their investigations.

In recent years, pure jump models have been widely used as an
alternative model for price process to the classical model, which
has a continuous martingale component; see Todorov and Tauchen
(\citeyear{TodTau10}) and references within.
% Pure-jump models are comprised solely of jumps.
The idea behind the pure-jump\vadjust{\goodbreak} modeling is that small jumps can
eliminate the need for a continuous martingale. The class of
pure-jump models is extremely wide. It includes the normal inverse
Gaussian [Rydberg (\citeyear{Ryd97}), Barndorff-Nielsen (\citeyear{Bar97,Bar98})], the
variance gamma [Madan, Carr and Chang (\citeyear{MadCarCha98})], the CGMY model of Carr
et al. (\citeyear{Caretal02}), the time-changed Levy models of Carr et
al. (\citeyear{Caretal03}), the COGARCH model of Kl{\"u}ppelberg, Lindner and Maller (\citeyear{KluLinMal04})
for the financial prices, as well as the non-Gaussian
Ornstein--Uhlenbeck-based models of Barndorff-Nielsen and Shephard
(\citeyear{BarShe01}) and the L\'{e}vy-driven continuous-time moving average
(CARMA) models of Brockwell (\citeyear{Bro01}) for the stochastic volatility.
Pure-jump models have been extensively considered and used for
general options pricing [Huang and Wu (\citeyear{HuaWu04}), Broadie and Detemple
(\citeyear{BroDet04}), Levendorskii (\citeyear{Lev04}), Schoutens (\citeyear{Sch06}), Ivanov (\citeyear{Iva07})], and for
foreign exchange options pricing [Huang and Hung (\citeyear{HuaHun05}), Daal and
Madan (\citeyear{DaaMad05}), Carr and Wu (\citeyear{CarWu07})]. Other applications of pure-jump
models include reliability theory [Drosen (\citeyear{Dro86})], insurance
valuation [Ballotta (\citeyear{Bal05})] and financial equilibrium analysis
[Madan (\citeyear{Mad05})].

Given the wide usage of pure jump models, % has been established,
a natural question is: \textit{is there any statistical evidence from
the high-frequency financial data to support using the purely
discontinuous models alone without any continuous diffusion
components?}
% Put it another way: \textit{is the diffusion term necessary, given that
%the jumps are present?}
% The reason for the interest of the problem is several-fold,
% There are several other reasons why we are interested in the above
%test,
% from both theoretical and practical considerations:
%
The question is of significance from both theoretical and practical
viewpoints: % for the following reasons:
\begin{itemize}
\item
Many empirical evidences indicate that pure jump models
% (such as CGMY models)
can fit the data well; see, for example, Cont and Tankov (\citeyear{ConTan04}), and Carr
et al. (\citeyear{Caretal02}) and references therein. Therefore, it would be
of theoretical interest to establish some statistical tests for this
purpose.

% find out whether the underlying asset returns are truly from a pure
% jump model or a mixture model containing both jumps and a diffusion.

\item
Given the existence of jumps, pure jump models are typically easier
to handle than mixture models in practice, and a preferred choice to
mixture models for users.
% For example, a $\beta$-stable process contains only one
% parameter, which can be estimated by the MLE.
% the CGMY model contains only four parameters and there are many
% available methods to estimate them.
% Under these circumstances, pure jump models would be a preferred
% choice to mixture models for users.
However, before using a pure jump model, one must check its
validity.

%has to take of the diffusion term as well as the jumps, and the
%mixing of both components make the model estimation a lot more
%difficult. ??????
%
%given the jumps are present. }

\item
Various jump models have been well studied in the literature, as
mentioned earlier. Should we decide to use pure jump models, we
would have an array of available tools at our disposal.

\item
Many results are strongly model dependent, and any model
mis-specification could have a severe effect on the results.
Therefore, it is imperative to choose the best possible model, and
model selection is very critical.
\end{itemize}

To put our question into a mathematical context, suppose that the
price process~$Y$ is a jump diffusion process of the form
%
%e1.1 ###
\begin{equation}
Y_t = X_t+J_t, \label{model1}
\end{equation}
for $t\in[0,T]$ with $X_t$ and $J_t$ being the continuous and
discontinuous (or jump)
components, defined % respectively
as
%
%e1.2 ###
\begin{eqnarray}\label{semmtg}
X_t &=& Y_0 + \int^t_0 b(X_s)\,ds + \int^t_0 \sigma(X_s) \,dW_s,
\nonumber
\\[-10pt]
\\[-10pt]
\nonumber
J_t &=& \int^t_0\int_{|x|\le1} x(\mu-\nu)(ds,dx) +
\int^t_0\int_{|x|>1}x\mu(ds,dx),
\end{eqnarray}
where $b$ and $\sigma$ are some deterministic functions such that
$X$ has unique weak solution, $\mu$ is the jump measure, with $\nu$
its predictable compensator; for details on jump diffusion
processes, see Jacod and Shiryaev (\citeyear{JacShi03}).
% All these are defined on a stochastic basis $(\Omega,\mathcal{F},$ $ {
%
Under this framework, the above question is tantamount to testing\vspace*{-1pt}
%
%e1.4 ###
%e1.3 ###
\begin{eqnarray}
H_0\dvtx \quad \int_0^T \sigma^2(X_s) \,ds &>& 0, \qquad  \mbox{(i.e.,
diffusion effect is present),} \label{testing.problem0}\hspace*{-35pt}
\\
H_1\dvtx\quad  \int_0^T \sigma^2(X_s) \,ds &=& 0, \qquad  \mbox{(i.e.,
diffusion effect is not present)},
\label{testing.problem}\hspace*{-35pt}
\end{eqnarray}
given the jump component $J_t$ is present. Note that, under $H_0$,
$Y_t$ is a mixture model of diffusion and jumps while, under $H_1$,
it is a pure jump model.

Cont and Mancini (\citeyear{ConMan}) and A\"{\i}t-Sahalia and Jacod (\citeyear{AtSJac10})
considered the above test using threshold power variation. They
assumed a general continuous semi-martingale form of $X$ as opposed
to a diffusion form in the present paper. However, to perform their
test, one needs to impose the condition that $J$ is of finite
variation [e.g., Theorem 2 of A\"{\i}t-Sahalia and Jacod (\citeyear{AtSJac10})].
This restriction rules out some interesting models used in finance,
where the jumps are shown to be of infinite variation, as done in
A\"{\i}t-Sahalia and Jacod (\citeyear{AtSJac09N1}), Zhao and Wu (\citeyear{ZhaWu09}) and some other
references mentioned earlier.

%serve as local alternatives since the small jumps behave like that
%of a continuous martingale driven by Brownian motion more than jump
%processes with finite variation. On the other hand, the test based
%on threshold power variation has less power as jumps behave less
%actively, see Table 5 in the simulation studies. This gives the
%impression that when the pure jump process is further away from the
%local alternatives, the less power the test has. Simulation also
%show that the test proposed in A\"{\i}t-Sahalia and Jacod (2010)
%perform not well in power for finite sample.
%}
%%

In this paper, we propose a simple-to-use, general purpose and yet
powerful goodness-of-fit test for differentiating a pure jump model
from a mixture model. The CLTs are also derived for the test
statistics under $H_0$, regardless whether the jump component is of
finite or infinite variation. In that aspect, our proposed test
works more generally than those proposed earlier by Cont and Mancini
(\citeyear{ConMan}) and A\"{\i}t-Sahalia and Jacod (\citeyear{AtSJac10}). Even for the
situations where tests by Cont and Mancini (\citeyear{ConMan}) and
A\"{\i}t-Sahalia and Jacod (\citeyear{AtSJac10}) are applicable, our numerical
results also show the superior performance of our proposed test.

The paper is organized as follows. In Section~\ref{methodology}, we
give some motivations via a simple example and then formally
introduce our test statistics. Asymptotic results are derived in
Section~\ref{mainresults}. Some review of alternative tests are
given in
Section~\ref{comparison}. Numerical studies are given % relegated to
in Section~\ref{simulation0}.
% , where comparisons are made with some alternative tests.
A real example is studied in Section~\ref{real}. Some discussion on
microstructure noise is given in Section~\ref{conclusion}. All
technical proofs are postponed in the \hyperref[proof]{Appendix}.

%f1 ###
\begin{figure}

\includegraphics{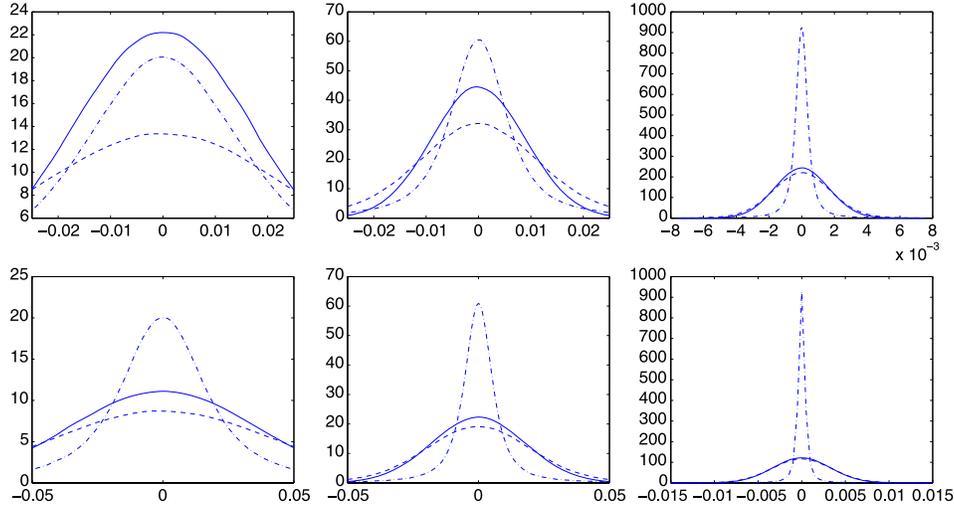}

\caption{Smoothed histograms for the increment of the mixture model
(- -), pure jump model (-$\cdot$), and diffusion term alone (-). From left
to right, the sample sizes are 195, 780 and 23,400, respectively.
From top to bottom, $\sigma=$ 0.25 and 0.5, respectively.}
\label{fig1}
\end{figure}

% \section{The setup} % problem}

Throughout the paper, the available data set is denoted as
$\{Y_{t_i}; 0 \le i\le n\}$ in the fixed interval $[0, T]$, which
is discretely sampled from $Y$. For simplicity, we assume that
$\{Y_{t_i}; 0 \le i\le n\}$ are equally spaced in $[0, T]$, that is,
$t_i=i\Delta_n$ with $\Delta_n=T/n$ for $0 \le i\le n$. Denote the
$j$th one-step increment by
\[
\Delta^n_j Y=Y_{t_j}-Y_{t_{j-1}}, \qquad 1 \le i\le n.
\]

%s2 ###
\section{Test statistics} \label{methodology}

We start with a simple motivating example first and then introduce
our test statistics for testing
(\ref{testing.problem0}) and (\ref{testing.problem}).

%s2.1 ###
\subsection{A simple motivating example} \label{simpleexample}

We draw two respective samples $\{Y_{t_i};\allowbreak  0 \le i\le n\}$ from
the following two models:
\begin{eqnarray*}
 H_0\dvtx\quad   Y_t &= &\sigma W_t + S_t^\beta \qquad  (\textit{a mixture model}), \\
 H_1\dvtx\quad   Y_t &=& S_t^\beta\qquad   (\textit{a pure jump model}),
\end{eqnarray*}
where $W_t$ and $S_t^\beta$ are a standard Brownian motion and a
symmetric $\beta$-stable L\'{e}vy process, respectively. So the
mixture model contains an extra continuous component $\sigma W_t$,
in comparison with the pure jump model. For illustration, we take
$T=1$, $\beta=1.25$ and $\sigma=0.25$, $0.5$.

The smoothed histograms (done by $10^6$ replications) of the
increments $\{\Delta_j^n Y, 1\le j\le n \}$ under the two models are
plotted in Figure~\ref{fig1} for sample sizes $n=195$, $780$, and
$23\mbox{,}400$, which corresponds to sampling every 2 minutes, 30 seconds,
and every second in a 6.5 hour trading day. From Figure~\ref{fig1},
we can see some very clear patterns:

\begin{enumerate}[(1)]
\item[(1)]
For small sample size $n$ and small $\sigma$, it is difficult to
distinguish the~models under $H_0$ and $H_1$ (the dashed line and
dash-dotted line). However, as $n$ and/or $\sigma$ increases, the
difference is more significant under~$H_0$ and~$H_1$.\looseness=-1\vadjust{\goodbreak}

\item[(2)]
The differences between the normal histogram and the mixture one
(the solid line and the dashed line) are small in all cases and
become even more negligible as $n$ increases. Literally, the jump
component has been ``absorbed'' by the diffusion component in the
center.

\item[(3)]
For fixed $\sigma$, as the sample size $n$ increases, the
differences between models under $H_0$ and $H_1$ are getting
sharper. Take $n=23\mbox{,}400$ and $\sigma=0.5$, for example. The histogram
under $H_1$ (dash-dotted line) shows a very narrow peak around the
origin, while the histogram under $H_0$ (the dashed line) stays
rather flat.
\end{enumerate}

The example shows that there is a huge difference around the origin
between the models under $H_0$ and $H_1$. If we use the number of
``small'' increments as an indicator, $U_n=\sum_{i=1}^n
\{|\Delta_i^nY|\le u_n \}$ for some $u_n$, then it relies heavily
on whether the diffusion is present or not, particularly when the
sample size $n$ gets large. To give a better idea, some values of
$U_n$ under the above two models are presented in Table~\ref{T1} when
$n=23\mbox{,}400$. The drastic difference for the two models strongly
suggests that we might be able to use~$U_n$ to test whether the
diffusion is present or not.

%t1 ###
\begin{table}
\caption{Numbers of increments $\le\alpha\Delta_n^{\varpi}$ for $Y$,
$W$ and $S^{\beta}$, where $\alpha=2$, $\varpi=1$ and $\Delta
_n=1/23\mbox{,}400$. The numbers
are averaged over 500 replications} \label{T1}
\begin{tabular*}{\textwidth}{@{\extracolsep{\fill}}lccc@{}}
\hline
$\mathbf{Parameter}$ & $\bolds{\#\{|\Delta_i^nY|\leq\alpha\Delta_n^{\varpi}\}}$ &
$\bolds{\#\{|\Delta_i^nW|\leq\alpha\Delta_n^{\varpi}\}} $ & $\bolds{ \#\{|\Delta
_i^nS|\leq\alpha\Delta_n^{\varpi}\}}$\\
\hline
$\beta=1.50$  &   $408$   &   $488$   &   $\phantom{16\mbox{,}}942$    \\
$\beta=1.00$  &   $485$  &  $489$  &  $16\mbox{,}491$ \\
$\beta=0.50$  &  $487$  &  $487$ &  $23\mbox{,}313$ \\
\hline
\end{tabular*}
\end{table}

%s2.2 ###
\subsection{Test statistics} \label{test_stat}

Let us return to the testing problem given in
(\ref{testing.problem0}) and (\ref{testing.problem}). We observe from
Section~\ref{simpleexample} that the increments from a pure jump
model and a mixture model have fundamentally different behavior
around the centers of their distributions. Namely, the distribution
for the increments from a pure jump model shows a much higher peak
in the center than that from a mixture model. In other words, the
number of small increments from a pure jump model is far greater
than than that from a mixture model. This suggests that we might use
the number of small increments
\[
U(\Delta_n) =:  U(\alpha, \Delta_n, \varpi,
T)=\sum_{i=1}^{[T/\Delta_n]}I(|\Delta^n_iY|\le
\alpha\Delta_n^{\varpi}),
\]
to define a test statistic. Note that $U(\Delta_n)$ simply counts
the number of increments smaller than $\alpha\Delta_n^{\varpi}$,
where $\alpha>0$ and $\varpi>1/2$. (Here, we suppress the dependence
on $\alpha$, $\varpi$ and $T$ for convenience.)

Under some mild conditions (given in Section~\ref{mainresults}), the
behaviors of $U(\Delta_n)$ are different under $H_0$ and $H_1$. Here
is a heuristic argument. Under $H_0$, we have $\Delta^n_iY \approx
\sigma(X_{t_{i-1}}) \Delta_i^nW$, and hence
\begin{eqnarray*}
EU(\Delta_n) &\approx&
% \sum^{[T/\Delta_n]}_{i=1}P(|\Delta^n_iY|\leq\alpha\Delta_n^{\varpi})
\sum^{[T/\Delta_n]}_{i=1} E P_{t_{i-1}} \bigl( |\Delta_i^nW| \leq
\alpha\Delta_n^{\varpi}/\sigma(X_{t_{i-1}}) \bigr)\\
&\approx&
2\alpha\phi(0)\Delta_n^{\varpi-3/2} T \int_0^T E\sigma^{-1}(X_s) \,ds,
\end{eqnarray*}
where $\phi(x)$ is the density of the standard normal r.v.,
and $P_{t_{i-1}}$ is the probability conditioned at time $t_{i-1}$.
Consequently, we have $U(\Delta_n)$ is of order
$\Delta_n^{-3/2+\varpi}$ under $H_0$. Similarly, we can show that
$U(\Delta_n)$ is of order $\Delta_n^{-(1+1/\beta)+\varpi}$ under
$H_1$. Clearly, we have $\Delta_n^{-3/2+\varpi} \ll
\Delta_n^{-(1+1/\beta)+\varpi}$. That is, there are far more small
increments under the pure jump model ($H_1$) than those under the
mixture model ($H_0$), which agrees well with the above motivating
example. Then, we will reject $H_0$ (a mixture model) in favor of
$H_1$ (a pure jump model), if $U(\Delta_n)$ is large enough.

From both Proposition~\ref{propo} and (\ref{prelln1}) in the
\hyperref[proof]{Appendix}, we see that the probability limit of $U(\Delta_n)$ depends
on unknown population quantities, and hence can not be directly used
for our testing purposes. To get around the problem, we adopt the
same strategy as in Zhang, Mykland and A\"{\i}t-Sahalia (\citeyear{ZhaMykAtS05}) and
A\"{\i}t-Sahalia and Jacod (\citeyear{AtSJac10}) by using a two-time scale test
statistic,
\[
V_n:=\frac{U(\Delta_n)}{U(k\Delta_n)},
\]
where $ U(k\Delta_n) \,{=:}\, U(\alpha, k\Delta_n, \varpi, T) \,{=}\,
\sum_{i=1}^{[T/(k\Delta_n)]} I(|\Delta^n_{(i-1)k+1}Y \,{+}\,
\cdots\,{+}\,
\Delta^n_{ik}Y|\,{\le}\allowbreak\alpha(k\Delta_n)^{\varpi})$. As can be seen
from (\ref{lln}) below, the distribution of $V_n$ is model-free
under $H_0$, and we can reject $H_0$ (a mixture model) in favor of
$H_1$ (a~pure jump model) if $V_n > C$ for some critical value
$C>0$.

We end the section by pointing out some differences between the
above test and the one by A\"{\i}t-Sahalia and Jacod (\citeyear{AtSJac10}). The
test statistic given in~(16) and (19) of A\"{\i}t-Sahalia and Jacod
(\citeyear{AtSJac10}) is based on the truncated $p$th power variations while our
test statistic given by $U(\Delta_n)$ and $V_n$ is simply based on
the number of small increments. Further comparisons will be made
later in the paper.

% Sections~\ref{comparison}-\ref{real}.

%s3 ###
\section{Main results} \label{mainresults}

We first list some assumptions and then present the main results.

%s3.1 ###
\subsection{Model assumptions} \label{setup}

Recall that $Y_t = X_t + J_t$. Assume that $Y$ is defined on a
filtered probability space $(\Omega, \mathcal{F}^Y,
\mathcal{F}^Y_t)$, where $\mathcal{F}^Y_t$ is the history of $Y$ up
to time $t$.\vadjust{\goodbreak}

\begin{assumption} \label{assump1}
$J_t$ has a jump measure $\mu(dx, dt)$ with compensator $\nu(\omega,
dx, dt)=dtF_t(\omega, dx)$, such that, for all $(\omega, t)$, we
have $F_t=F_t' + F_t''$, where:
\begin{enumerate}[(1)]
\item[(1)]
$F_t'$ has the form
\[
F_t'(dx)=\frac{1+|x|^{\gamma}f( x)}{|x|^{1+\beta}} \bigl[a^{(+)}I(0<x\le
\varepsilon^{+})+a^{(-)}I(-\varepsilon^{-}\le x< 0) \bigr] \,dx,
\]
for some positive constants $a^{(+)}$, $a^{(-)}$, $\gamma$,
$\varepsilon^{+}$ and $\varepsilon^{-}$ and some bounded function $f(x)$,
satisfying $1+|x|^{\gamma}f(x)> 0,  |f(x)|\le L.$
\item[(2)]
$F_t''$ is a singular measure with respect to $F_t'$, satisfying $
\int_R(|x|^{\beta^{\prime}}\wedge1)F^{\prime\prime}_t(\omega,\allowbreak
dx)\le L.$
\end{enumerate}
\end{assumption}

\begin{assumption} \label{assump2}
$X$ and $J$ are mutually independent.
\end{assumption}

\begin{assumption} \label{assump3}
$b(\cdot)$ is a bounded continuous functions, $\sigma(\cdot)$ is
bounded away from zero and infinity if it does not vanish and
$\sigma^{\prime}(\cdot)$ exists and is bounded.
\end{assumption}

Assumption~\ref{assump1} implies that the small jumps of $J$ form a
L\'{e}vy process with a $\beta$-stable-like L\'{e}vy density, while
almost no condition is placed on the large jumps of $J$, and
$F_t^{\prime\prime}$ could even be random. Assumption~\ref{assump1}
includes a rich class of models, like the variance gamma model, CGMY
model, tempered stable process, etc. Assumptions
\ref{assump2} and~\ref{assump3} are technical conditions.
%

%s3.2 ###
\subsection{Asymptotic results}

Let $\mathcal{N}(0,1)$ denote a standard Gaussian random variable. We
will use the stable convergence in law below, which is slightly
stronger than weak convergence; see, for example, Jacod and Shiryaev
(\citeyear{JacShi03}).

\begin{theorem} \label{Theorem1}
Suppose that $\varpi>\beta-1/2$ and that Assumptions 1--3 hold.
\begin{enumerate}[(1)]
\item[(1)]
We have
%
%e3.1 ###
\begin{equation} \label{lln}\qquad
V_n \rightarrow^P \cases{
k^{3/2-\varpi}, &  \quad $\mbox{under $H_0$},$\vspace*{2pt}\cr
k^{1+(1/\beta-\varpi)\wedge0}, &\quad $\mbox{under $H_1$ and
Assumption~\ref{assump4} below}.$}
\end{equation}
\item[(2)]
Let $k=2$. Under $H_0$, we have
\[
\Delta_n^{({\varpi-3/2})/{2}} (V_n-k^{3/2-\varpi}) \longrightarrow
\sigma\mathcal{N}(0,1) \qquad \mbox{stably},
\]
where $\mathcal{N}(0,1)$ % $z$ is a standard Gaussian variable
is independent of $Y$ and % $\mathcal{F}^Y$ and
\[
\sigma^2=\frac{(1+k^{3/2-\varpi})k^{3-2\varpi}}{2\alpha\phi(0)\int^T_0
\sigma^{-1}(X_s) \,ds}.
\]
\end{enumerate}
\end{theorem}

%%
%% The condition $k=2$ is not needed in (\ref{lln}).
%The stable convergence is slightly stronger than weak convergence;
%see, e.g., Jacod and Shiryaev (2003).}
%%
%}

To apply Theorem~\ref{Theorem1}, one needs to estimate the unknown
$\sigma^2$. However, in view of Proposition~\ref{propo} in the
\hyperref[proof]{Appendix} and the stable convergence, we have the following.\vadjust{\goodbreak}
% corollary.

\begin{corollary} \label{coro1}
Assuming the same assumptions as in Theorem~\ref{Theorem1}, we have
\[
\Delta_n^{({\varpi-3/2})/{2}} (V_n-k^{3/2-\varpi})/\widehat\sigma
\longrightarrow_d \mathcal{N}(0,1)\qquad  \mbox{under $H_0$},
\]
where
\[
\widehat\sigma^2 = \frac{(1+k^{3/2-\varpi})
k^{3-2\varpi}}{\Delta_n^{3/2-\varpi} U(\Delta_n)}.
\]
\end{corollary}

From Corollary~\ref{coro1}, at significance level $\theta$, we can
reject $H_0$ if $ V_n > k^{3/2-\varpi} + z_{1-\theta}
\Delta_n^{3/4-\varpi/2} \widehat\sigma$ and $P(\mathcal{N}(0,1) >
z_{1-\theta}) = \theta$.
It follows from Corollary~\ref{coro1} that the size of the above
test
% in (\ref{aaaa})
is asymptotically $\theta$.

A slight variant of the test statistic $V_n$ can be given below. Let
\[
\tilde V_n = \frac{U(\Delta_n)}{U_L(2\Delta_n)},
\]
where $U_{ L}(2\Delta_n)=[U(2\Delta_n)+U^{\prime}(2\Delta_n)]/2$ and
\begin{eqnarray*}
U^{\prime}(2\Delta_n) &=&
\sum^{[T/(2\Delta_n)]-1}_{i=1}I\bigl(|\Delta^n_{2i+1}Y+\Delta^n_{2i}Y|\leq
\alpha(2\Delta_n)^{\varpi}\bigr), \\
U(2\Delta_n) &=&
\sum^{[T/(2\Delta_n)]}_{i=1}I\bigl(|\Delta^n_{2i}Y+\Delta^n_{2i-1}Y|\leq
\alpha(2\Delta_n)^{\varpi}\bigr).
\end{eqnarray*}
In other words, we use linear combinations of $U(2\Delta_n)$ with
different starting time points instead of a single $U(2\Delta_n)$
starting from time $t_0$ when nonoverlapping two-step increments of
$Y$ are sampled. Similarly to Corollary~\ref{coro1}, we can easily derive the
following result.

\begin{corollary} \label{cor2}
Assuming the same assumptions as in Theorem~\ref{Theorem1}, we have
\[
\Delta_n^{({\varpi-3/2})/{2}} (\tilde V_n-2^{3/2-\varpi})/\tilde
\sigma\longrightarrow_d \mathcal{N}(0,1) \qquad \mbox{under $H_0$},
\]
where
\[
\tilde\sigma^2 = \frac{ U(\Delta_n)+2^{3/2-\varpi}U_L(2\Delta_n)/2
}{ \Delta_n^{3/2-\varpi}
U_L(2\Delta_n)^2}.
\]
\end{corollary}

Our final decision rule is: at significance level $\theta$, we
reject $H_0$ if
%
%e3.2 ###
\begin{equation}
\tilde V_n > \tilde C, \label{aaaaa}
\end{equation}
where $\tilde C=2^{3/2-\varpi} + z_{1-\theta}
\Delta_n^{3/4-\varpi/2} \tilde\sigma$ and $P(\mathcal{N}(0,1) >
z_{1-\theta}) = \theta$.
It follows from Corollary~\ref{cor2} that the size of the above test
in (\ref{aaaaa}) is asymptotically $\theta$.

\begin{remark}\label{rem1}
%In Theorem~\ref{Theorem1} and Corollaries~\ref{coro1}-\ref{cor2}, we
%require $\varpi$ to satisfy $\varpi>\beta-1/2$. Since $\beta\in[0,
%2)$, one can choose $\varpi=1.5$. This is of course a conservative
%choice, and will be adopted in the simulation studies later.
%}
The requirement $\varpi>\beta-1/2$ in Theorem~\ref{Theorem1} and
Corollaries~\ref{coro1}--\ref{cor2} can be easily satisfied by
choosing $\varpi=3/2$ as $\beta\in(0,2)$.
%
%we only require $\varpi>\beta-1/2$, and no upper bound is requested
%for $\varpi$, so conservatively, we can choose $\varpi=1.5$ since
%$\beta$ is at most 2. }
%
Moreover, whatever the value of $\varpi$, $H_0$ and $H_1$ can be
differentiated since $1+1/\beta>3/2$ for all $\beta\in(0, 2)$.

On the other hand, the behaviors of test statistics $S_n$ under
$H_0$ and~$H_1$ in A\"{\i}t-Sahalia and Jacod (\citeyear{AtSJac10}) depend on the
choice of $p$, that is, $2>p>1 \vee\beta$; see Theorem 1 in that
paper. A\"{\i}t-Sahalia and Jacod (\citeyear{AtSJac10}).
% to guarantee the consistency of $S_n$, the test statistic, under
%$H_1$,
% $p$ has to satisfy $2>p>1 \vee\beta$.
Since $\beta$ is unknown, it is difficult to choose $p$.
% can not be chosen to be 2
To be on the safe side, one might try to choose $\beta$ close to~2.
However, this will render the test with very low power since $S_n$
converges in probability to roughly the same limit 1 under~$H_0$ and~$H_1$.
\end{remark}

\begin{remark}\label{rem2}
In Theorem~\ref{Theorem1} and Corollaries
\ref{coro1}--\ref{cor2}, we have $\beta\in(0,2)$, and no further
restriction on $\beta$ is imposed, so that the jump component could
be of finite variation or infinite variation. By contrast, The CLT
under $H_0$ was developed by A\"{\i}t-Sahalia and Jacod (\citeyear{AtSJac10}) only
when $\beta<1$, namely when $J$ is of finite variation.
% Further discussion will be given in Section~\ref{comparison}.
%be controlled approximately by the central limit theorem given in
%Corollary~\ref{cor2}. }
\end{remark}

%s3.3 ###
\subsection{Asymptotic power}
Before discussing the power of our test statistic, we list one more
condition, which % Assumption~\ref{assump4}
basically assumes that the drift term is zero when $\beta\le1$. It
is a standard assumption
in the literature; see Jacod (\citeyear{Jac08}) and Woerner (\citeyear{Woe03}), and the
references therein. %?????? ??????? ??????? ??????.

% ?????? $\beta\in(0,2)$, the assumption sounds like that $\beta\le
%1$, anything wrong? ???????
%
\begin{assumption} \label{assump4}
% \item
If $\beta<1$, we assume that $b(\cdot)\equiv0$, and $\int_{|x|\leq
1}xF^{\prime}(dx)\equiv0$.
% \item
If $\beta=1$, we assume that $b(\cdot)\equiv0$ and $F^{\prime}(dx)$
is symmetric about 0.
% \end{enumerate}
\end{assumption}

%Assumption~\ref{assump4} basically assumes that the drift
%term is zero when $\beta\le1$. It is a standard assumption in the
%literature; see ?????? ??????? ??????? ??????. It could be relaxed
%to a great extent, e.g., in the simulation study we intentionally
%keep a nonzero drift, i.e., $b(\cdot)\not\equiv0$ for all $\beta
%assumed here in evaluating the power; see the proof of Theorem 2,
%which
%}

The next theorem gives the asymptotic power of our proposed test
(\ref{aaaaa}).

\begin{theorem} \label{th2}
Under Assumptions~\ref{assump1} and~\ref{assump4}, with prescribed
level $\theta$ and for $\varpi>1$, we have
\[
P( \tilde V_n > \tilde C | H_1 ) \longrightarrow1,
\]
that is, the asymptotic power is 1.
\end{theorem}

\begin{remark} \label{remark_power}
We end this section with some remarks on finite sample performance
of our test statistics.
% Theorem~\ref{th2} ensures that we will eventually be able to
%differentiate between the mixture model and
% pure jump model. The question is how the test performs in finite
%sample situations.
Intuitively, the closer $\beta$ gets to 2, the more the pure jump
process behaves like a diffusion process; thus, the more difficult
it is to tell their difference apart, the less power our test will
have. Similarly, the closer~$\beta$ gets to 0, the more power our
test will have. In fact, simple algebra yields
$\tilde{C}-2^{3/2-\varpi}=O_p(\Delta_n^{(1+1/\beta-\varpi)/2})$,
from which we can see that, as $\beta$ becomes closer to $0$, the
power of our test increases soon. This is further confirmed in our
simulation studies given later.
\end{remark}

%s4 ###
\section{A review of other approaches} \label{comparison}

The testing problem considered in this paper has also been
considered earlier by % others, including
Cont and Mancini (\citeyear{ConMan})\vadjust{\goodbreak} and A\"{\i}t-Sahalia and Jacod (\citeyear{AtSJac10}). Since
the work in both papers is similar, we will only review the test by
A\"{\i}t-Sahalia and Jacod (\citeyear{AtSJac10}) (hereafter AJ's test) below.

The building block of the AJ's test is based on the truncated
$p$-power variation,
%
%e4.1 ###
\begin{equation}
B(p, u_n,
\Delta_n)=\sum^{[t/\Delta_n]}_{i=1}|\Delta^n_iY|^pI(|\Delta_i^nY|\le
u_n), \label{tpv}
\end{equation}
where $p \in(1,2)$, and $u_n$ satisfies % is a threshold level
%satisfying
$u_n/\Delta_n^{{\rho}_{-}} \rightarrow0$,
$u_n/\Delta_n^{{\rho}_{+}} \rightarrow\infty$, for some $0 \le
{\rho}_{-} < {\rho}_{+} < 1/2$.
Similarly to Zhang, Mykland and A\"{\i}t-Sahalia (\citeyear{ZhaMykAtS05}), A\"{\i}t-Sahalia and Jacod
(\citeyear{AtSJac10}) defined a two-time scale estimator
\[
S_n=\frac{B(p, u_n, \Delta_n)}{B(p, u_n, k\Delta_n)}\qquad
\mbox{for an integer $k\ge2$},
\]
and showed that
%
%e4.2 ###
\begin{equation}
S_n\rightarrow^P \cases{
k^{1-p/2}, &  \quad$\mbox{under $H_0$,}$ \vspace*{2pt}\cr
1, &  \quad $\mbox{under $H_1$,  if $2>p>1\vee\beta$ and ${\rho}_{+} \le
{(p-1)}/{p}$}$} \label{llntpv1}\hspace*{-35pt}
\end{equation}
and that when $\beta<1$,
%
%e4.3 ###
\begin{eqnarray}
(S_n - k^{1-p/2})/{\sqrt{v_n}} \longrightarrow_d \mathcal{N}(0,1)
 \qquad\mbox{under $H_0$,} \label{snkp}
\end{eqnarray}
where $v_n^2=CB(2p, u_n, \Delta_n)/B(p, u_n, \Delta_n)^2$ for some
constant $C$. Noting\break $k^{1-p/2} > 1$, one would reject $H_0$ if $S_n
\le C_0$, for some $C_0$ determined from the CLT. A\"{\i}t-Sahalia
and Jacod (\citeyear{AtSJac10}) also showed that the asymptotic power of this test
is 1.

We make several remarks regarding the AJ's test:
% In our view, there are several shortcomings with the A\"{\i}t-Sahalia
%and
% Jacod's test, despite its nice theoretical properties given above.
%
\begin{itemize}
\item
From (\ref{llntpv1}),
% the choice of $p$ depends on the value of unknown $\beta$, which
% can be particularly difficult if $\beta$ is close to $2$.
the behaviors of test statistics $S_n$ under $H_0$ and $H_1$ depend
on the choice of $p$, that is, $2>p>1 \vee\beta$.
% to guarantee the consistency of $S_n$, the test statistic, under
%$H_1$,
% $p$ has to satisfy $2>p>1 \vee\beta$.
Since $\beta$ is unknown, it is difficult to choose~$p$.
% can not be chosen to be 2
To be on the safe side, one might try to choose $\beta$ close to 2.
However, this will render the test with very low power since $S_n$
converges in probability to roughly the same limit 1 under $H_0$ and
$H_1$.

\item
The CLT under $H_0$, (\ref{snkp}), was established in
A\"{\i}t-Sahalia and Jacod (\citeyear{AtSJac10}) only for the case $\beta<1$,
namely when $J$ is of finite variation. However, when $\beta>1$,
that is, when $J$ is of infinite variation, no CLT is available, and
hence the size of the test cannot be controlled for that case. This
rules out some interesting applications when $\beta>1$.

\item
For $\beta\in(0,1)$, where the CLT is available for AJ's test, we
might expect that AJ's test should have very good power,
particularly as $\beta$ gets smaller toward 0.
% This is based on the similar argument to that in Remark
% should, in principle, apply to the test proposed in A\"{\i}t-Sahalia
%and Jacod (2010).
However, our simulation studies give some counterintuitive
results; see Table~\ref{table6}. %It is not clear exactly why
%this happens; further investigation into this should be
%worthwhile.

% This forms a contrary to the power performance of the test in
% A\"{\i}t-Sahalia and Jacod (2010) mentioned earlier in the
% Introduction.

\end{itemize}

%}

%s5 ###
\section{Numerical studies} \label{simulation0}

In this section, we conduct simulations to evaluate the performance
of our proposed test statistics, and make some comparisons with that
of A\"{\i}t-Sahilia and Jocod (\citeyear{AtSJac10}).\vadjust{\goodbreak}
% We also apply the test to some real data sets.

The test statistics $V_n$ and ${\widetilde V}_n$ involve choosing
the threshold level
$u_n=\alpha\Delta_n^{\varpi}$. %, which can be carried out as follows.
In view of the requirement $\varpi>\beta-1/2$, a conservative choice
of $\varpi$ would be 1.5. To compensate for the conservative choice
of $\varpi$, we choose a relatively large $\alpha$ by
$\alpha_n=\delta(\log n)^{\kappa}$ for some positive constants
$\delta$ and $\kappa$. This choice will not affect any of the
asymptotic results in the paper.

%But in this case few data would be retained since $u_n$ might be
%quite small even if $\alpha$ is chosen
%relatively large. %, say $20$-$50$.
%To remedy this, we suggest to replace $\alpha$ by
%$\alpha_n=\delta(\log n)^{\kappa}$ for some positive constants
%$\delta$ and $\kappa$, which will not effect any of the asymptotic
%results in the paper.
%}

% \subsection{Simulation studies} \label{simulation1}

Assume that the data generating process under the null and
alternative hypotheses are, respectively,
%
%e5.2 ###
%e5.1 ###
\begin{eqnarray}
H_0\dvtx \quad   Y_t &=& X_t + \theta^{\prime} S_{\beta,t}, \label{sim_h0}
\\
H_1\dvtx \quad  Y_t&=&\exp(-\gamma t) + 0.5 S_{\beta, t},
\label{sim_h1}
\end{eqnarray}
where $X_t$ is an Ornstein--Urlenbeck process. $dX_t=-X_t \,dt + dW_t$,
% $$ dX_t=-\gamma X_tdt+\sigma dW_t, $$
% where $\gamma=1$, $\sigma=1$
and $W$ is a standard Brownian motion, and $S_{\beta}$ is a
symmetric $\beta$-stable process. Let $T=1$, $\theta^{\prime}=0.5$.
Also we take $n=1560,  2340,  4680, 11\mbox{,}700  \mbox{ and }  23\mbox{,}400$,
corresponding to an intra day data set recorded every 15, 10, 5, 2
and 1 seconds in a 6.5-hour trading day, respectively. We will
simulate $10\mbox{,}000$ samples from each model above.

\subsection*{Asymptotic sizes}

% First, we consider the asymptotic size.
Fix the nominal level $\theta=5\%$, so the critical value is
$z_{0.95}=1.645$. The size of the test is calculated by the
percentage of samples such that~(\ref{aaaaa}) holds true over
$10\mbox{,}000$ samples.

%t2 ###
\begin{table}
\caption{Sizes of the test $(\%)$ under different $n$'s and $\beta$'s,
 $(\delta=2$, $\kappa=2$)} \label{T3}
\begin{tabular*}{\textwidth}{@{\extracolsep{\fill}}lcccccccc@{}}
\hline
\textbf{Value of} $\bolds{\beta}$ & $\bolds{1.2}$ &
$\bolds{1.3}$ & $\bolds{1.4}$ & $\bolds{1.5}$ & $\bolds{1.6}$ & $\bolds{1.7}$ & $\bolds{1.8}$ & $\bolds{1.9}$\\
\hline
$n=1560$ & $2.97$ & $4.18$ & $4.71$ & $4.08$ & $4.39$
& $3.99$ & $4.13$ & $4.27$\\
$n=2340$ & $3.94$ &
$4.42$ & $4.41$ & $4.27$ & $4.55$ & $4.42$ & $4.25$ & $4.37$\\
$n=4680$ & $3.98$ &
$4.36$ & $4.46$ & $4.82$ & $5.31$ & $4.81$ & $4.93$ & $4.56$\\
$n=11\mbox{,}700$ & $4.34$ &
$4.50$ & $4.94$ & $5.06$ & $4.93$ & $5.06$ & $4.78$ & $4.29$\\
\hline
\end{tabular*}
\end{table}

Table~\ref{T3} reports the asymptotic sizes for different sample
sizes. From the table, we see that the type I error is well
controlled by $5\%$; as the sample size $n$ increases, the
asymptotic sizes become closer to the true size $5\%$.

Table~\ref{T4} reports the asymptotic sizes across different
threshold levels which reflect the number of effective data. It
shows that control of type I error is not affected much by changes
of $\delta$.

\subsection*{Asymptotic power}

We also consider the power performance of ${\widetilde V}_n$. The
power of the test is the percentage of samples with (\ref{aaaaa})
violated over $10\mbox{,} 000$ samples. The results are listed in Table
\ref{T5} for different values of $\beta$.

From Table~\ref{T5}, it is clear that,
% We make the following remarks.
% from Table~\ref{T5}:
% \begin{enumerate}
% \item
as the sample size $n$ increases, the test becomes more powerful
overall, as expected.
%
% \item
The test is powerful especially when $\beta$ is away\vadjust{\goodbreak} from $2$. When
$\beta$ approaches $2$, the power gradually diminishes. This is
easily understandable as in this case the behavior of the
discontinuous process resembles that of a Brownian motion. This can
also be seen from (\ref{lln}).

% \end{enumerate}

%t3 ###
\begin{table}
\caption{Sizes of the test $(\%)$ under different $\delta$'s and $\beta
$'s $(n=23\mbox{,}400$, $\kappa=2)$}\label{T4}
\begin{tabular*}{\textwidth}{@{\extracolsep{\fill}}lcccccccc@{}}
\hline
\textbf{Value of} $\bolds{\beta}$ & $\bolds{1.2}$ &
$\bolds{1.3}$ & $\bolds{1.4}$ & $\bolds{1.5}$ & $\bolds{1.6}$ & $\bolds{1.7}$ & $\bolds{1.8}$ & $\bolds{1.9}$\\
\hline
 $\delta=1.0$ & $4.64$ &
$4.22$ & $4.21$ & $4.68$ & $4.61$ & $4.16$ & $3.98$ & $4.03$\\
$\delta=1.5$ & $4.23$ &
$4.16$ & $4.65$ & $4.54$ & $4.62$ & $4.97$ & $4.35$ & $4.65$\\
$\delta=2.0$ & $3.94$ &
$4.42$ & $4.41$ & $4.27$ & $4.55$ & $4.42$ & $4.25$ & $4.37$\\
$\delta=2.5$ & $4.30$ &
$4.62$ & $4.16$ & $4.38$ & $4.47$ & $4.49$ & $4.65$ & $4.08$\\
\hline
\end{tabular*}
\end{table}

%t4 ###
\begin{table}
\caption{Powers of the test $(\%)$ under different $n$'s and $\beta
$'s $(\delta=2$, $\kappa=2$)}\label{T5}
\begin{tabular*}{\textwidth}{@{\extracolsep{\fill}}lcccccccc@{}}
\hline
\textbf{Value  of}  $\bolds{\beta}$ & $\bolds{1.2}$ &
$\bolds{1.3}$ & $\bolds{1.4}$ & $\bolds{1.5}$ & $\bolds{1.6}$ & $\bolds{1.7}$ & $\bolds{1.8}$ & $\bolds{1.9}$\\
\hline
$n=1560$ & $100$ &
$100$ & $100$ & $94.10$ & $48.79$ & $19.86$ & \phantom{0}$9.85$ & $5.38$\\
$n=2340$ & $100$ &
$100$ & $100$ & $91.36$ & $46.46$ & $20.53$ & $10.43$ & $5.82$\\
$ n=4680 $ & $100$ &
$100$ & $100$ & $89.55$ & $48.52$ & $21.85$ & $17.55$ & $6.14$\\
$n=11\mbox{,}700$ & $100$ &
$100$ & $100$ & $92.78$ & $55.33$ & $26.54$ & $12.67$ & $7.13$\\
$n=23\mbox{,}400$ & $100$ &
$100$ & $100$ & $96.01$ & $63.80$ & $31.05$ & $14.59$ & $7.34$\\
\hline
\end{tabular*}
\end{table}

\begin{figure}[b]

\includegraphics{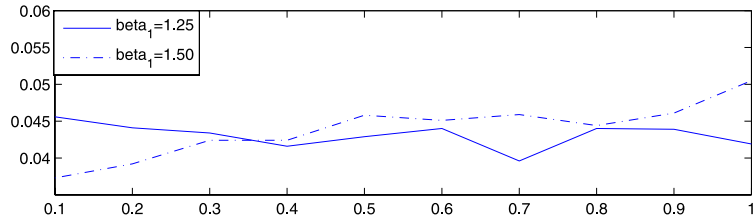}

\caption{Sensitivity plot of asymptotic sizes to choices of
$\theta^{\prime}$.} \label{fig4+}
\end{figure}

Finally, we examine the asymptotic sizes over different choices of
$\theta^{\prime}$. We fix $n=2340$, $\delta=2$, $\kappa=2$ and
$\theta=5\%$. In Figure~\ref{fig4+}, the asymptotic sizes for
$\beta=1.25$ and $1.5$ are plotted against $\theta^{\prime}$.
Clearly, the asymptotic sizes are not sensitive to choices of
$\theta^{\prime}$.

%f2 ###

\subsection*{Comparisons with AJ's test}

Now we compare the performance of our estimator $\tilde V_n$ with
that of AJ's estimator $S_n$, under the same settings as in~(\ref{sim_h0}) and (\ref{sim_h1}). However, since AJ's test is only
shown to be valid for the case $\beta\in(0,1)$ (i.e., the jump
process is of finite variation), our comparisons are also restricted
to that case. Tables~\ref{table5} and~\ref{table6} report the sizes
and powers of our test and AJ's test for various values of $\beta
\in(0,1)$, respectively.

\begin{itemize}
\item
For both tests, all sizes are close to to the nominal level, $5\%$,
with AJ's test being slightly closer overall.

\item
Our test outperforms AJ's in terms of power throughout. In fact, our
test has full power for all $\beta\in(0,1]$, even for sample size
$n=1560$. On the other hand, AJ's test has very low power in
detecting the alternatives for $\beta\le0.7$, even when $n=23\mbox{,}400$.
\end{itemize}

The very low powers of AJ's test for small $\beta$ came as a
surprise to us. Some more detailed analysis suggests that the reason
might be due to the large variation of $S_n$ for finite sample size
$n$ under $H_1$. More precisely, from~(69) in A\"{\i}t-Sahalia and
Jacod (\citeyear{AtSJac10}), we have $S_n=O_P(u_n^{\beta/2})$ under $H_1$.
%where $u_n$ is usually taken as of order $\Delta_n^{1/4}$.
So for finite sample $n$, $S_n$ may not be close to 0 for small
$\beta$, which often places the test statistic $S_n$ wrongly within
the acceptance region, resulting in low power. It also explains why
the problem is mostly pronounced if $\beta$ is closer to 0.

\begin{figure}

\includegraphics{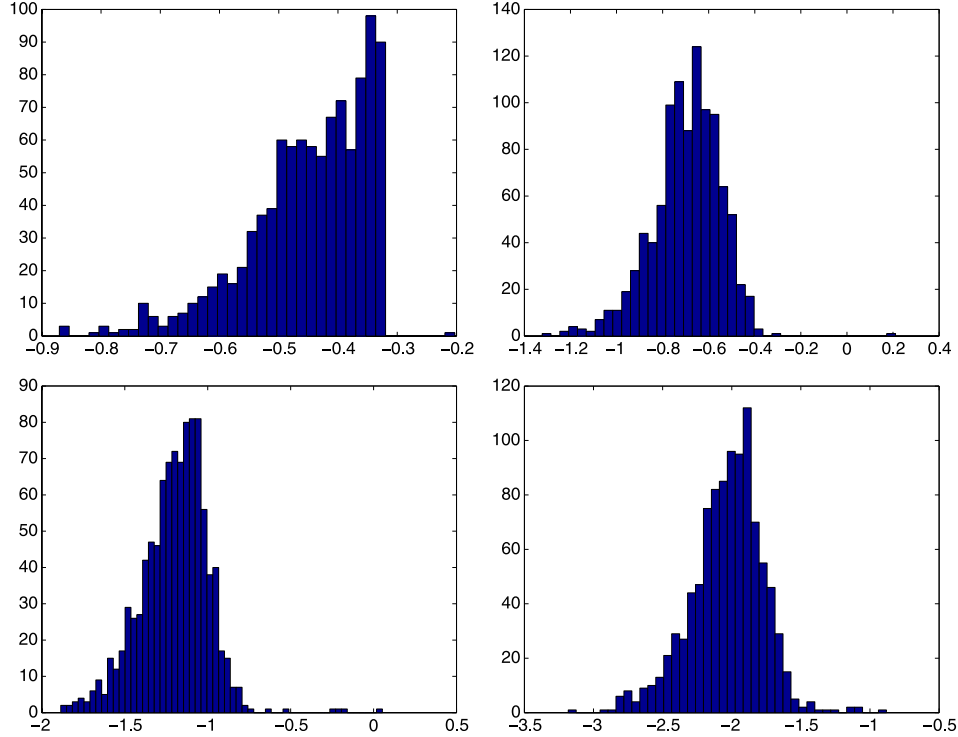}

\caption{Upper left panel: $\beta=0.25$; upper right panel:
$\beta=0.50$; lower left panel: $\beta=0.75$; lower right
panel: $\beta=1.00$.} \label{fig2}
\end{figure}

Figure~\ref{fig2} displays the histograms of the studentized $S_n$
as given in (\ref{snkp}) when $n=4680$. We see that values of
studentized $S_n$'s are seldom less than $z_{0.05}=-1.645$, except
in the case $\beta=1$.

%t5 ###
\begin{table}
\caption{Size of our test v.s. that of AJ's test, $(\%)$,  $\delta
=2$,  $\kappa=2$}\label{table5}
\begin{tabular*}{\textwidth}{@{\extracolsep{\fill}}lccccccccccc@{}}
\hline
& $\bolds{\beta}$ & $\bolds{0.1}$ & $\bolds{0.2}$ & $\bolds{0.3}$ & $\bolds{0.4}$ & $\bolds{0.5}$ & $\bolds{0.6}$ &
$\bolds{0.7}$ & $\bolds{0.8}$ & $\bolds{0.9}$ & $\bolds{1.0}$ \\
\hline
$n=1560$ & $\tilde V_n$ & $4.19$ & $3.81$ & $4.18$ & $3.90$ & $4.08$ &
$3.84$ & $3.62$ & $3.59$ &$4.06$ & $4.11$\\
& $AJ$ & $4.08$ & $4.23$ & $4.22$ & $3.94$ & $3.99$ & $4.29$ &
$4.24$ & $4.10$ & $4.37$ & $4.41$\\
$n=2340$ & $\tilde V_n$ & $4.31$ & $4.37$ & $3.97$ & $3.98$ & $4.16$ &
$4.10$ & $4.18$ & $3.93$ & $4.11$ & $4.32$\\
 & $AJ$ & $4.36$ & $4.33$ & $4.23$ & $3.89$ & $4.67$ & $4.47$ &
$4.61$ & $4.29$ & $4.54$ & $4.47$\\
$n=4680$ & $\tilde V_n$ & $4.07$ & $4.38$ & $4.01$ & $4.57$ & $3.87$ &
$4.24$ & $4.34$ & $4.18$ & $4.52$ & $4.43$\\
 & $AJ$ & $4.52$ & $4.53$ & $4.33$ & $4.64$ & $4.56$ & $4.43$ &
$4.41$ & $4.71$ & $4.89$ & $ 4.76$\\
\hline
\end{tabular*}
\end{table}
%

%t6 ###
\begin{table}[b]
\tabcolsep=0pt
\caption{Power comparisons of our test v.s. that of AJ's test when
$\beta\in(0,1]$, $(\%)$, $\delta=2$, $\kappa=2$} \label{table6}
\begin{tabular*}{\textwidth}{@{\extracolsep{\fill}}lcd{3.0}d{3.0}d{3.0}d{3.2}d{3.0}d{3.2}d{3.2}d{3.2}d{3.2}d{3.2}@{}}
\hline
& $\bolds{\beta}$ & \multicolumn{1}{c}{$\bolds{0.1}$} & \multicolumn{1}{c}{$\bolds{0.2}$} & \multicolumn{1}{c}{$\bolds{0.3}$} & \multicolumn{1}{c}{$\bolds{0.4}$} &
\multicolumn{1}{c}{$\bolds{0.5}$} & \multicolumn{1}{c}{$\bolds{0.6}$} & \multicolumn{1}{c}{$\bolds{0.7}$} &
\multicolumn{1}{c}{$\bolds{0.8}$} & \multicolumn{1}{c}{$\bolds{0.9}$} & \multicolumn{1}{c@{}}{$\bolds{1.0}$}\\
\hline
$n=1560$ & $\tilde V_n$
& 100 & 100 & 100 & 100 & 100 & 100 & 100 & 100 & 100 & 100 \\
& $AJ$ &  0  &  0  &  0  &  0  &  0  &  0.01  &  0.25  &  2.03
&  7.61  &  24.03 \\
$n=2340$ & $\tilde V_n$ &  100  &  100  &  100  &  100  &  100  &  100
&  100  &  100  &  100  &  100 \\
& $AJ$ &  0  &  0  &  0  &  0.01  &  0  &  0.07  &  0.43  &
 3.32  &  14.82  &  48.36 \\
$n=4680$ & $\tilde V_n$ &  100  &  100  &  100  &  100  &  100  &  100
&  100  &  100  &  100  &  100 \\
& $AJ$ &  0  &  0  &  0  &  0  &  0  &  0.08  &  0.97  &  8.66
&  48.26  &  97.02 \\
$n=23\mbox{,}400$ &  $\tilde V_n$  &  100  &  100  &  100  &  100  &  100  &
 100  &  100  &  100  &  100  &  100 \\
& $AJ$ &  0  &  0  &  0  &  0  &  0  &  0.30  &  6.40  &  85  &
 100  &  100 \\
\hline
\end{tabular*}
\end{table}

%
%f3 ###

\subsection*{Sensitivity to model misspecification of our test}
\begin{figure}

\includegraphics{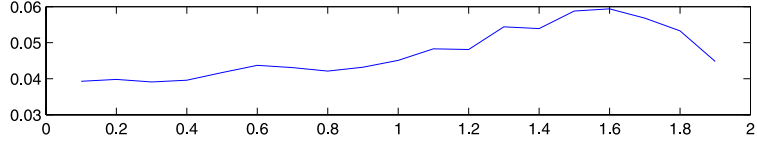}

\caption{The data generating process is the combination of a
stochastic volatility process and a standard symmetric stable
process.} \label{fig3}
\end{figure}

In the model assumptions, we assumed a local volatility function.
Now we conduct a simulation study to check the sensitivity of our
test to model misspecification; see Figure~\ref{fig3}. Instead of
using an Ornstein--Urlenbeck process as the continuous part of the
full model, we use a stochastic volatility process here, that is,
$dX_t=\sigma_t \,dW_t$ with $\sigma_t=v_t^{1/2}$,
$dv_t=\kappa(\eta-v_t)\,dt+\gamma v_t^{1/2}\,dB_t$, $E[dW_t\,dB_t]=\rho
\,dt$. We take $\eta=1/16$, $\gamma=0.5$, $\kappa=5$, $\rho=-0.5$. We
use $\theta^{\prime}S_{\beta, t}$ as the jump process as in last two
simulations. Now we fix $n=23\mbox{,}400$, $\delta=1$, $\kappa=2$,
$\theta^{\prime}=0.25$ and $\theta=5\%$. All simulations are run
$10\mbox{,}000$ times. From Figure~\ref{fig3}, the asymptotic sizes are not
much affected by using a stochastic volatility model as the
continuous part.

%f4 ###

%s6 ###
\section{A real data set analysis} \label{real}

In this section, we implement our test to some real data sets. We
use the stock price records of Microsoft (MFST) in there trading
days, Nov. 1, Dec. 1 and Dec. 11 in the year 2000. All data sets
are from the TAQ database. For prices recorded simultaneously, we
use their averages. To weaken the possible effect from
microstructure noise, we sparsely sample observations every~10
seconds and the sample sizes for the aforementioned three trading
days are $1343$, $1701$ and $1253$, respectively. Finally, we take the
logarithm of the sparsely sampled prices and use the log prices to
calculate the test statistics. We set $T=1$ (day) consisting of 6.5
hours of trading time.

\begin{figure}[b]

\includegraphics{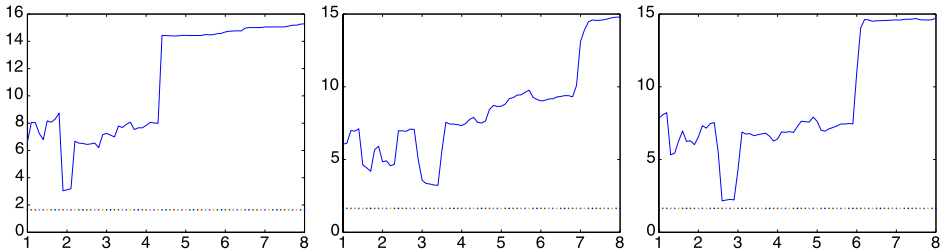}

\caption{The statistics evaluated over different values of $\delta$.
From left to right: test statistics for 01, Nov., 01, Dec. and 11,
Dec., respectively. The horizontal axis stands for the value of
$\delta$ while the vertical axis stands for the value of the test
statistics.}\label{fig4}
\end{figure}

We now discuss how to choose the parameters $\delta$, $\kappa$ and
$\varpi$. As argued theoretically at the beginning of Section
\ref{simulation0}, we fix $\varpi=1.5$. Since $\kappa$ and $\delta$
are dependent parameters, we fix $\kappa=2$ and consider a grid of
points of $\delta$ such that
%
%e6.1 ###
\begin{equation}
\delta(\log{n})^{\kappa}\Delta_n^{\varpi}\leq
\hat{\sigma}^*\Delta_n^{1/2}, \label{delta}
\end{equation}
where $\hat{\sigma}^*$ is approximately the averaged standard
deviation of the diffusion component of one 10-second log return in
case the diffusion term exists in the underlying dynamics.
Mathematically, $\hat{\sigma}^*$ is defined as
\[
\hat{\sigma}^{*2}=:\frac{1}{T}\sum(\Delta^n_iX)^2I(|\Delta^n_iX|\leq
\Delta_n^{1/4})\rightarrow^P \frac{1}{T}\int^T_0\sigma^2(X_s)\,ds;
\]
see Jacod (\citeyear{Jac08}) for example. In virtue of (\ref{delta}), we can
choose $\delta$ conservatively as the grid points from 1 to 8 with
equal step length 0.1 for all three data sets. The plots are
displayed in Figure~\ref{fig4}.

%f5 ###

From the plots, the observed test statistics are all larger than
$1.645$. Therefore we can reject the existence of the diffusion
component and simply use a~pure jump model to characterize the
underlying dynamics of the prices for those three days.

%s7 ###
\section{Discussions on microstructure noise} \label{conclusion}

%In this paper, we studied the issue of whether asset price
%processes can be modeled by pure jump processes instead of
%jump-diffusion models. A very simple statistical test was devised
%for this purpose.
%% for pure discontinuity of stochastic processes against jump
%% diffusion models using high frequency data.
%Compared with earlier tests, our test works for general jump
%processes with finite or infinite variation. It can control the
%type I error probability accurately and has asymptotic power 1.
%% it is also powerful.
%Simulations justify the good finite-sample performance of the test
%statistics.
%}

It is widely accepted nowadays that microstructure noise is present.
Various methods have been studied to handle the issue of the
microstructure noise in the context of the integrated volatility
estimation for high-frequency data. See, for example,
A\"{\i}t-Sahalia, Mykland and Zhang (\citeyear{AitMykZha05}), Zhang, Mykland and A{\"{\i}}t-Sahalia (\citeyear{ZhaMykAtS05}),
Zhang (\citeyear{Zha06}), Fan and Wang (\citeyear{FanWan07}), Podolskij and Vetter (\citeyear{PodVet09})
% Barndorff-Nielsen \textit{et al.} (2009a) (2009b)
and Jacod et al. (\citeyear{Jacetal09}), among others. A very effective
technique in handling microstructure noise is the so-called
``pre-averaging method''; see Jacod et al. (\citeyear{Jacetal09}) and Podolskij
and Vetter (\citeyear{PodVet09}).

Suppose that the observation at time $t_i$ is
\[
Z_{t_i}=Y_{t_i} + \varepsilon_{t_i},\qquad   i=1,\ldots,n,
\]
where $Y_t$ is an unobserved % underlying
semi-martingale of the form (\ref{model1}) and (\ref{semmtg}), and~$\varepsilon_{t_i}$ with
mean 0 and variance $\sigma^2$ %\sim N(0, \sigma^2)$
is the microstructure noise at time~$t_i$. We wish to test
(\ref{testing.problem0}) and (\ref{testing.problem}), that is,
whether $Y_t$ can be modeled as a~pure jump
process, or not.

%%
%H_0: \ Y_t=W_t+S_t,  \mbox{vs}   H_1: \ Y_t=S_t,
%%
%where $W$ and $S$ are respectively the standard Brownian motion and
%symmetric $\beta$-stable process, and $\varepsilon$, $W$, and $S$ are
%independent.
%
%}

So far, we have not seen any work in the testing framework in the
presence of microstructure noise.
% To find out how the ``pre-averaging method" work in this context,
We now apply the simplest pre-averaging technique as follows. We
first separate the full data set $Z_{t_i}$, $1\leq i\leq n$ into
$n/M$ nonoverlapped blocks,
\[
\{ Z_{t_1}, \ldots, Z_{t_M} \},  \ldots,  \bigl\{ Z_{t_{kM+1}},
\ldots, Z_{(k+1)M} \bigr\}.
\]
Then within each block, we take the average of all $K$-step
increments, that is,
\begin{eqnarray}
\overline{Z}_j=\frac{1}{n/M-K+1}\sum^{(k+1)M}_{i=kM+K+1}(Z_{t_{i}}-Z_{t_{i-K}})
:=\overline{X}_j+\overline{J}_j+\overline{\varepsilon}_j,
\nonumber
\\[-8pt]
\eqntext{j=1,\ldots,n/M.}
\end{eqnarray}
Simple calculation yields $\overline{X}_j=O_p(M^{-1/2}),  \overline{J}_j=O_p(M^{-1/\beta}),
\overline{\varepsilon}=\break O_p((M/n)^{1/2}).$ By properly tuning $M$, for example,
$M=o(n^{1/2})$, one could make $\overline{\varepsilon}_j$
asymptotically negligible. Based on the modified data set
$\overline{Z}_1, \ldots, \overline{Z}_M$, the test statistics can be
defined (similarly to $V_n$) as
\[
\overline{V}_n=\frac{\overline{U}(\Delta_M)}{\overline{U}(k\Delta_M)},
\]
where $\Delta_M=T/M$, $\overline{U}(\Delta_M)$ and
$\overline{U}(k\Delta_M)$ are defined as $U(\Delta_M)$ and
$U(k\Delta_M)$ by replacing $Y_i$ with $\overline{Z}_i$ and by
replacing $\Delta_n$ by $\Delta_M$, for example, $U(\Delta_M) =
\sum^M_{i=1}I(|\overline{Z}_i|\leq\alpha(\Delta_M)^{\varpi})$.

%%
%U(\Delta_M) &=& \sum^M_{i=1}I(|\overline{Z}_i|\leq\alpha(\Delta
%_M)^{\varpi}), \\
%U(2\Delta_M) &=&
%U^{\prime}(2\Delta_M) &=&
%U_L(2\Delta_M) &=& (U(2\Delta_M)+U^{\prime}(2\Delta_M))/2.
%%
%}

Under appropriate conditions, the results obtained in the paper
should be expected to hold here as well, for instance,
% we expect to have
% As might be expected, we should have
%
%e7.1 ###
\begin{equation} \label{consistency}
\overline{V}_n \rightarrow^P \cases{
k^{1.5-\varpi}, &\quad
$\mbox{under  }H_0,$\vspace*{2pt}\cr
k^{1+(1/\beta-\varpi)\wedge0}, &\quad $\mbox{under } H_1$.}
\end{equation}

\begin{figure}

\includegraphics{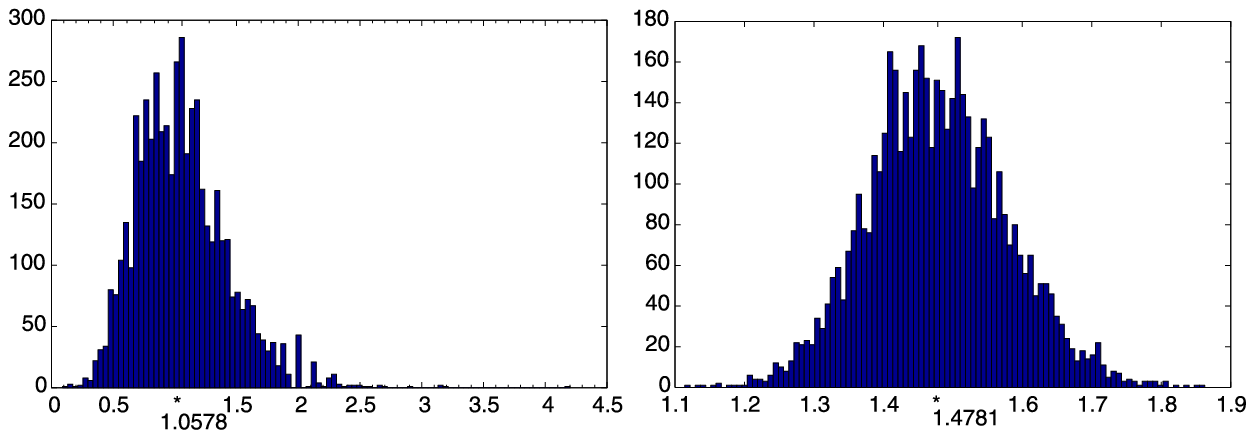}

\caption{Histograms of $\overline{V}_n$.}\label{microfig}
\end{figure}

Let us conduct a simple simulation study to investigate the
feasibility of the test statistic $\overline{V}_n$.
% we conduct a simple simulation.
Take $Y_t=W_t+S_t$ under $H_0$ and $Y_t=S_t$ under $H_1$, where
$W_t$ and~$S_t$ are
% We generate the data set by simulating
a standard Brownian motion and a symmetric Cauchy process (i.e.,
$\beta=1$), respectively. Also take $\sigma^2 \sim N(0, \sigma^2)$
with $\sigma=0.01$. We let $T=1$, $n=23\mbox{,}400$ and $k=2$. We further
take $M=234$, $K=50$, $\alpha=9$, $\varpi=1.5$. Note that the choice
of $M=234$ corresponds to taking averages about every 4 minutes. The
simulation is repeated 5000 times. Each time, we calculate
$\overline{V}_n$ both under $H_0$ and~$H_1$. Their histograms under
$H_0$ and $H_1$ are plotted in Figure~\ref{microfig}.
% The first figure is the histogram of $\tilde{V}_n$ under $H_0$
% while the second under $H_1$.
% Figure~\ref{microfig1}.

From Figure~\ref{microfig}, we see that the means of
$\overline{V}_n$ under $H_0$ and $H_1$ (marked by $*$ in the
horizontal axis) are 1.0578 and 1.4781, respectively. These are
rather close to the asymptotic values 1 and 1.414, given by
(\ref{consistency}). Note that the effective sample size after
pre-averaging is $23\mbox{,}400/234 = 100$, a rather small sample size for
this testing purpose. This explains partly why the variances of
histograms plots are rather large, and there are substantial overlaps
between the plots under $H_0$ and $H_1$. If we choose $M=120$, or
even $60$, then the histograms under $H_0$ and $H_1$ will become
thinner and more easily separable.

% The following figures show that (\ref{consistency}) is
% desirable. The means are marked by $*$ in the horizontal axis.
%
%f6 ###

The above simple simulation study suggests that the pre-averaging
method would work well in handling microstructure noise in the
testing problems. Of course, there remain many theoretical and
practical issues to be resolved. For example, we need to establish a
CLT under $H_0$; to study its asymptotic power; to find a
data-driven method to determine parameter $M$, etc. We will pursue
these and other related issues in our future work.

%%
%%
%%
%%
%asymptotically 1.
%%
%%
%good in finite sample. This is not the case when $M$ is not very
%large, see for example Figure 1.
%%
%}

\begin{appendix}
\section*{Appendix} \label{proof}

In the sequel, $C$ will denote a constant which may take different
values in different places, and $\chi$ is an arbitrarily small
positive number. Also, $P_{t_{i-1}}$ and $E_{t_{i-1}}$ denote
probability and expectation given time $t_{i-1}$, respectively.

%s7.1 ###
\subsection{\texorpdfstring{Proof of Theorem \protect\ref{Theorem1}}{Proof of Theorem 1}}

Let $\sigma_0^2=2\alpha\phi(0)\int^T_0 \sigma^{-1}(X_s) \,ds$. Now
\begin{eqnarray*}
&& \Delta_n^{({\varpi-3/2})/{2}} ( V_n-k^{3/2-\varpi} ) \\
&&\qquad=
\frac{\Delta_n^{({\varpi-3/2})/{2}}(\tilde{U}(\Delta_n)-\sigma_0^2)
-k^{({3/2-\varpi})/{2}}(k\Delta_n)^{({\varpi-3/2})/{2}}
(\tilde{U}(k\Delta_n)-\sigma_0^2)}{k^{\varpi-3/2}\tilde{U}(k\Delta_n)}\\
&&\qquad:=\frac AB.
\end{eqnarray*}
By Proposition~\ref{propo} (below), $A\rightarrow^S \sigma_0
z_1-k^{{(3/2-\varpi)}/{2}}\sigma_0 z_2$ with $z_1$ and $z_2$
independent Gaussian random variables independent of
$\mathcal{F}^Y$, while $B\rightarrow^P k^{\varpi-3/2}\sigma_0^2$,
which is random but depending only on $\mathcal{F}^Y$. Then Theorem
\ref{Theorem1} is proved.

Now, we prove Proposition~\ref{propo}, in which we need the following two
lemmas. Lemma~\ref{le2} implies that the proportion of paths of a jump
diffusion process having ``small'' increments is the same as that of
the diffusion component. This has its own interest.

\begin{lemma}\label{le1} Let $A_i = \{\omega\dvtx  |\Delta^n_iX+x|\le\alpha\Delta
_n^{\varpi}
\}$.
% Under the conditions in Theorem 1,
%
\begin{enumerate}[(1)]
\item[(1)] For $|x|<\Delta_n^{1/2}$,
$ |P_{t_{i-1}}(A_i) -
\frac{2\alpha\phi(0)\Delta_n^{\varpi-1/2}}{\sigma(-x+X_{t_{i-1}})}|\le
C\Delta_n^{\varpi}( x^2\Delta_n^{-3/2}+\Delta_n^{-\chi}
). $
\item[(2)] For any $x\in R/\{0\}$, we have
$ P_{t_{i-1}}(A_i)\le C\Delta_n^{\varpi-1/2}. $
\end{enumerate}
\end{lemma}
\begin{pf}
Define $f(x)=\int^x_0 \sigma^{-1}(y)\, dy$, then
$f'(x)=\sigma^{-1}(x)$ and $f''(x)=-{\sigma'(x)}/{\sigma^2(x)}$. So
$f(x)$ is strictly increasing.
Let $\Xi_t=f(X_t)$, or equivalently,
% by the monotonicity of $f(x)$ we have
$X_t=f^{-1}(\Xi_t)$. By It\^{o}'s formula,
% \begin{equation}
% df(X_t)=(\frac{b(X_t)}{\sigma(X_t)}-\frac{1}{2}\sigma^{
% \label{timechange}
% \end{equation}
% Then (\ref{timechange}) can be rewritten as
%
\setcounter{equation}{0}
%e7.2 ###
\begin{eqnarray}\label{timechanged diffusion}
d\Xi_t&=&\biggl(\frac{b\circ f^{-1}(\Xi_t)}{\sigma\circ
f^{-1}(\Xi_t)}-\frac{1}{2}\sigma^{\prime}\circ f^{-1}(\Xi_t) \biggr)
\,dt + dW_t
\nonumber
\\[-8pt]
\\[-8pt]
\nonumber
&:= &\overline{b}\circ f^{-1}(\Xi_t)\,dt+dW_t.
\end{eqnarray}
%
% Since $X_t$ is a weak solution to (\ref{diffusion}) without
% explosion within $[0, T]$ and $\sigma(\cdot)$ is strictly positive,
% $\Xi$ is a weak solution to (\ref{timechanged diffusion}) without
% explosion within $[0, T]$, or equivalently
% Clearly, $\Xi$ is a weak solution to the martingale problem
% $MP(\overline{b}\circ f^{-1}, 1)$ under $(\Omega, \mathcal{F},
Let $\mathcal{F}^{\prime}_{t}=\mathcal{F}_{t+t_{i-1}}$ where
$t_{i-1}$ is the $(i-1)$th observation time defined at the end of the
\hyperref[sec1]{Introduction},
% then $\Xi_{t+t_{i-1}}$, $t\geq0$ is a solution to
% $MP(\overline{b}\circ f^{-1}, 1)$ under $(\Omega, \mathcal{F},
and $\tilde{W}_t=W_{t+t_{i-1}}$, $t\geq0$.
% On one hand, for any $A\in\mathcal{F}^{\prime}_s$ and any $B\in
% \begin{equation}
% E(\tilde{W}_t; A\cap B)=E(\tilde{W}_s; A\cap B), \ \ 0\le s\le t.
% \label{coditioning}
% \end{equation}
% which says that
It is easy to see that $\tilde{W}$ is a martingale under $(\Omega,
\mathcal{F}, \mathcal{F}^{\prime}_{t}, P_{t_{i-1}})$ with
% . On the other hand, the
quadratic variation $t$.
% of $\tilde{W}$ in $[0, t]$ is still $t$,
Thus by L\'{e}vy's characterization theorem, $\tilde{W}_t$ is a
Brownian motion under $(\Omega, \mathcal{F},
\mathcal{F}^{\prime}_{t}, P_{t_{i-1}})$.
By the Girsanov theorem,
% Take $b=-\overline{b}\circ f^{-1}$ in Theorem 5.1 of Durrett (2004),
% by the boundedness of the diffusion coefficient,
there exists a probability measure $Q_{t_{i-1}}$, locally equivalent
to $P_{t_{i-1}}$, satisfying
%
%e7.3 ###
\begin{equation}\qquad
\frac{dQ_{t_{i-1}}}{dP_{t_{i-1}}}\bigg|_{\mathcal{F}^{\prime}_t}
= \exp\biggl( -\int^t_0\overline{b}(X_{s+t_{i-1}})\,d\tilde{W}_s-\frac
{1}{2}\int^t_0\overline{b}^2(X_{s+t_{i-1}})\,ds \biggr),
\label{density}
\end{equation}
such that
% under $Q_{t_{i-1}}$, $\Xi_{t+t_{i-1}}$, $t\geq0$, is a solution to
%MP(0, 1). That is,
$\Xi_{t+t_{i-1}}$, $t\geq0$, is a Brownian motion under
$Q_{t_{i-1}}$.

Now
%
%e7.4 ###
\begin{eqnarray}\label{beforechange}\qquad
P_{t_{i-1}}(A_i)&=&P_{t_{i-1}}(-x-\alpha\Delta_n^{\varpi}+X_{t_{i-1}}\le
X_{t_i}\le-x+\alpha\Delta_n^{\varpi}+X_{t_{i-1}})
\nonumber\\
&=& P_{t_{i-1}}\bigl(f(-x+X_{t_{i-1}}-\alpha\Delta_n^{\varpi})
\nonumber
\\[-8pt]
\\[-8pt]
\nonumber
&&\phantom{P_{t_{i-1}}\bigl(}\le
f(X_{t_i})\le
f(-x+X_{t_{i-1}}+\alpha\Delta_n^{\varpi})\bigr)\\
&=& P_{t_{i-1}}(l_i\le\Xi_{t_{i}}-\Xi_{t_{i-1}}\le u_i),\nonumber
\end{eqnarray}
where $l_i=f(-x+X_{t_{i-1}}-\alpha\Delta_n^{\varpi})-f(X_{t_{i-1}})$
and $u_i=f(-x+X_{t_{i-1}}+\alpha\Delta_n^{\varpi})-f(X_{t_{i-1}})$.
Taking $t=t_{i}-t_{i-1}$ in (\ref{density}), we then have
%
%e7.5 ###
\begin{eqnarray}\label{afterchange}
P_{t_{i-1}}(A_i)&=& \int_{A_i}\frac{dP_{t_{i-1}}}{dQ_{t_{i-1}}}\,dQ_{t_{i-1}}
\nonumber\\
&= &\int_{A_i} \exp\biggl( \int^{t_{i}}_{t_{i-1}}\overline{b}\circ
f^{-1}(\Xi_s)\,dW_s\\
&&\phantom{\int_{A_i} \exp\biggl(}{}+\frac{1}{2}\int^{t_i}_{t_{i-1}}(\overline{b}\circ
f^{-1})^2(\Xi_s)\,ds \biggr)  \,dQ_{t_{i-1}}.\nonumber
\end{eqnarray}
Since $|\exp(x)-1|\le2|x|$ for $|x|\le\log2$, by boundedness of
the diffusion coefficient, L\'{e}vy's theorem of continuity modulus
and change of time,
\begin{eqnarray*}
&&\bigl| P_{t_{i-1}}(A_i) -
\bigl(\Phi\bigl(u_i/\sqrt{\Delta_n}\bigr)-\Phi\bigl(l_i/\sqrt{\Delta_n}\bigr)\bigr)\bigr|\\
&&\qquad\le
C\Delta_n^{1/2-\chi}\bigl(\Phi\bigl(u_i/\sqrt{\Delta_n}\bigr)-\Phi\bigl(l_i/\sqrt{\Delta_n}\bigr)\bigr),
\end{eqnarray*}
for any arbitrarily small $\chi>0$. On the other hand, by the mean
value theorem,
%
%e7.6 ###
\begin{equation}\qquad
\Phi\bigl(u_i/\sqrt{\Delta_n}\bigr)-\Phi\bigl(l_i/\sqrt{\Delta_n}\bigr)=\phi(\xi)
\frac{(u_i-l_i)}{\sqrt{\Delta_n}}=2\alpha\Delta_n^{\varpi-1/2}\frac{\phi
(\xi)}{\sigma(\eta)},
\label{concentration}
\end{equation}
where $\xi\in\frac{1}{\sqrt{\Delta_n}}(l_i, u_i)$ and
$\eta\in(-x+X_{t_{i-1}}-\alpha\Delta_n^{\varpi},
-x+X_{t_{i-1}}+\alpha\Delta_n^{\varpi})$. Then as
$\Delta_n\rightarrow0$, by Assumption~\ref{assump3}, we have
%
%e7.7 ###
\begin{equation}
|\sigma(\eta)- \sigma(-x+X_{t_{i-1}})|\le C\Delta_n^{\varpi}.
\label{sigma}
\end{equation}
Since as $n$ large enough $|u_i|\le C|x|$ and $|l_i|\le C|x|$ which
yields that $\xi\in\frac{1}{\sqrt{\Delta_n}}(-C|x|, C|x|)$. Then,
since $\phi^{\prime}(0)=0$ and $\phi^{\prime\prime}(\cdot)$ is
bounded,
%
%e7.8 ###
\begin{equation}
|\phi(\xi)-\phi(0)|\le C(x)^2\Delta_n^{-1}. \label{phi}
\end{equation}
The combination of (\ref{beforechange})--(\ref{phi}) completes the proof.
\end{pf}

\begin{lemma}\label{le2} Let $B_i = \{\omega\dvtx  |\Delta^n_i Y| \le\alpha\Delta
_n^{\varpi}
\}$. Then,
\[
\biggl|P_{t_{i-1}}(B_i) -
\frac{2\alpha\phi(0)}{\sigma(X_{t_{i-1}})}\Delta_n^{\varpi-1/2}\biggr|\le
C(\Delta_n^{\varpi-1/2+1-\beta/2}+\Delta_n^{\varpi-\chi}).
\]
\end{lemma}
\begin{pf} We write
%
%e7.9 ###
\begin{eqnarray}\label{p1p2}
&&P_{t_{i-1}}(B_i)\nonumber\hspace*{-35pt}\\
% &=& P_{t_{i-1}}(|\Delta^n_iY|\le
% \alpha\Delta_n^{\varpi})=P_{t_{i-1}}(|\Delta^n_iX+
% \alpha\Delta_n^{\varpi})\nonumber\\
&&\quad= \biggl(\int_{|x|< \sqrt{\Delta_n}}+\int_{|x|\geq
\sqrt{\Delta_n}}\biggr)P_{t_{i-1}}(|\Delta^n_iX+x|\le
\alpha\Delta_n^{\varpi})\,dP_{t_{i-1}}(\Delta^n_iJ\leq x
)
\hspace*{-35pt}\\
&&\qquad=:P_{i, 1}+P_{i, 2}.\nonumber\hspace*{-35pt}
\end{eqnarray}
Since $J$ is purely discontinuous, we can take the exponent in (64)
of A\"{\i}t-Sahalia and Jacod (\citeyear{AtSJac09N2}) as $1/2$, and then by Lemma~\ref{le1},
%
%e7.10 ###
\begin{equation}
P_{i, 2}\le
C\Delta_n^{\varpi-1/2}P_{t_{i-1}}\bigl(|\Delta^n_iJ|\geq\sqrt{\Delta
_n}\bigr)\le
C\Delta_n^{\varpi-1/2+1-\beta/2}. \label{p2}
\end{equation}
By Lemma~\ref{le1},
%
%e7.11 ###
\begin{equation}
P_{i, 1}=\int_{|x|<
\sqrt{\Delta_n}}\frac{2\alpha\phi(0)\Delta_n^{\varpi-1/2}}{\sigma
(-x+X_{t_{i-1}})}\,dP_{t_{i-1}}(\Delta^n_iJ\leq
x)+R_{n, i}. \label{p1}
\end{equation}
%
% For the same reason as in (\ref{p2}), we can take the exponent of
% the threshold as $1/2$ in Lemma 5 of Jing et al (2009), then
Similarly, we can obtain
%
%e7.12 ###
\begin{eqnarray}
|R_{n, i}| &\le& \int_{|x|<
\sqrt{\Delta_n}}C[(x)^2\Delta_n^{\varpi-3/2}+\Delta_n^{\varpi-\chi}]
\,dP_{t_{i-1}}(\Delta^n_iJ\leq x) \nonumber\\
&=&
C\Delta_n^{\varpi-3/2}E_{t_{i-1}}(\Delta^n_iJ)^2I\bigl(|\Delta^n_iJ|<
\sqrt{\Delta_n}\bigr)+C\Delta_n^{\varpi-\chi}\\
&\le&
C(\Delta_n^{\varpi-1/2+1-\beta/2}+\Delta_n^{\varpi-\chi}).\nonumber
\end{eqnarray}
Since $|x|< \sqrt{\Delta_n}$,
$|\sigma(-x+X_{t_{i-1}})-\sigma(X_{t_{i-1}})|\le C\sqrt{\Delta_n}$,
by boundedness of the diffusion coefficient,
%
%e7.13 ###
\begin{eqnarray}\label{prep3}
&&\int_{|x|<
\sqrt{\Delta_n}}2\alpha\phi(0)\Delta_n^{\varpi-1/2}\biggl(\frac{1}{\sigma
(-x+X_{t_{i-1}})}\!-\!\frac{1}{\sigma(X_{t_{i-1}})}\biggr)\,dP_{t_{i-1}}
(\Delta^n_iJ\!\leq\!x)\nonumber
\hspace*{-38pt}\\[-8pt]\\[-8pt]
\nonumber
&&\quad\le C\Delta_n^{\varpi}.\hspace*{-38pt}
\end{eqnarray}
On the other hand, as in (\ref{p2}), we have
%
%e7.14 ###
\begin{eqnarray}\label{p3}
&&\biggl|\frac{2\alpha\phi(0)\Delta_n^{\varpi-1/2}}{\sigma(X_{t_{i-1}})}
\biggl(\int_{|x|\le
\sqrt{\Delta_n}}\,dP_{t_{i-1}}(\Delta^n_iD\leq
x)-1\biggr)\biggr|
\nonumber
\\[-8pt]
\\[-8pt]
\nonumber
&&\qquad\le C\Delta_n^{\varpi-1/2+1-\beta/2}.
\end{eqnarray}
Combining (\ref{p1}), (\ref{prep3}) and (\ref{p3}) gives
\[
\biggl|P_{i,
1}-\frac{2\alpha\phi(0)\Delta_n^{\varpi-1/2}}{\sigma(X_{t_{i-1}})}\biggr|\le
C(\Delta_n^{\varpi-1/2+1-\beta/2}+\Delta_n^{\varpi-\chi}),
\]
which together with (\ref{p2})
completes the proof.
\end{pf}

Define $\tilde{U}(\Delta_n)=\Delta_n^{3/2-\varpi}U(\Delta_n)$, and
so $\tilde{U}(k\Delta_n)=(k\Delta_n)^{3/2-\varpi}U(k\Delta_n)$. Then
we have
% then we have the following proposition
%
\begin{proposition} \label{propo} We have % If $k=2$,
\begin{eqnarray*}
&&\Delta_n^{{(\varpi-3/2)}/{2}}\bigl(\tilde{U}(\Delta_n)-\sigma_0^2,
k^{{(\varpi-3/2)}/{2}} [ \tilde{U}(k\Delta_n)-\sigma_0^2 ]
\bigr)\\
&&\qquad\rightarrow\sigma_0 (z_1, z_2),\qquad
\mbox{$\mathcal{F}^Y$-stably,}
\end{eqnarray*}
where $z_1$ and $z_2$ are two independent Gaussian variables
independent of $\mathcal{F}^{Y}$.
\end{proposition}
\begin{pf} Without loss of generality, assume $k=2$. Denote $I_i =
I(|\Delta^n_iY|\le
\alpha\Delta_n^{\varpi})$. In view of Lemma~\ref{le2},
%
%e7.15 ###
\begin{equation}
\Biggl|\Delta_n^{3/2-\varpi}\sum_{i=1}^{[T/\Delta_n]} E_{t_{i-1}} I_i
%(|\Delta^n_iY|\le\alpha\Delta_n^{\varpi})
-\sigma_0^2\Biggr|\le
C(\Delta_n^{1-\beta/2}+\Delta_n^{1/2-\chi}). \label{bias}
\end{equation}
Since $\chi$ could be made arbitrarily small, and
$\varpi>\beta-1/2$, or equivalently,
$1-\beta/2>\frac{3/2-\varpi}{2}$,
%$$\Delta_n^{\frac{\varpi-3/2}{2}}(\tilde{U}(\Delta_n)-\sigma_0^2
%)$$
%
%e7.16 ###
\begin{equation}
\quad\qquad\Delta_n^{({\varpi-3/2})/{2}} \bigl(\tilde{U}(\Delta_n)-\sigma
_0^2\bigr)
=\Delta_n^{{3/2-\varpi}/{2}} \sum^{[T/\Delta_n]}_{i=1} (I_i -
E_{t_{i-1}} I_i) +o(1). \label{var}
\end{equation}
Now the summands in (\ref{var}) are centered martingale difference
sequences w.r.t. $\mathcal{F}_{t_{i-1}}$, $1\le i\le[T/\Delta_n]$.
In view of Lemma~\ref{le2}, and making use of (\ref{bias}) again,
\begin{eqnarray*}
\Delta_n^{3/2-\varpi}\sum^{[T/\Delta_n]}_{i=1} E_{t_{i-1}} (
I_i - E_{t_{i-1}} I_i )^2
& =&
\Delta_n^{3/2-\varpi}\sum^{[T/\Delta_n]}_{i=1} E_{t_{i-1}} I_i
+ o_P(1)\\
 &=& \sigma_0^2+o_P(1).
\end{eqnarray*}
Since the indicator function is bounded, the Linderberg condition
for the~martingale central limit theorem holds automatically. Then
by~(\ref{bias}) and~(\ref{var}),
%
%e7.17 ###
\begin{equation}
\Delta_n^{{(\varpi-3/2)}/{2}}\bigl(\tilde{U}(\Delta_n)-\sigma
_0^2\bigr)\rightarrow
\sigma_0 z_1 \label{clt1}
\end{equation}
$\mathcal{F}^{Y}$-stably if the following holds [c.f. Theorem IX 7.28
of Jacod and Shiryaev (\citeyear{JacShi03})]: for any bounded martingale
$N\in\mathcal{F}^{Y}$
%
%e7.18 ###
\begin{equation}
\Delta_n^{{(3/2-\varpi)}/{2}} \sum^{[T/\Delta_n]}_{i=1}
E_{t_{i-1}} ( \Delta^n_i N ) I_i
% ( |\Delta^n_i Y| \le\alpha\Delta_n^{\varpi} )
\rightarrow^P 0. \label{stablecon}
\end{equation}
Since $\mathcal{F}^{Y}=\mathcal{F}^{X}\vee\mathcal{F}^{J}$, it
suffices to show (\ref{stablecon}) with $N$ replaced by $X$ and
$N_1\in\mathcal{F}^{J}$, respectively, where $N_1$ is a bounded
martingale. By L\'{e}vy's theorem of continuity modulus,
(\ref{bias}) and $\varpi>1/2$,
%
%e7.19 ###
\begin{eqnarray}
&&\Delta_n^{{(3/2-\varpi)}/{2}} \sum^{[T/\Delta_n]}_{i=1}
E_{t_{i-1}} ( \Delta^n_i X ) I_i
\nonumber
\\[-10pt]
\\[-10pt]
\nonumber
&&\qquad \le
C\Delta_n^{{(3/2-\varpi)}/{2}+1/2-\chi}\sum^{[T/n]}_{i=1}
E_{t_{i-1}} I_i
% (|\Delta^n_iY|\le\alpha\Delta_n^{\varpi})
\rightarrow^P 0.
\end{eqnarray}
Next, by independence of $X$ from $\mathcal{F^{J}}$ and Lemma~\ref{le1},
%
%e7.20 ###
\begin{eqnarray}
&&\Delta_n^{{(3/2-\varpi)}/{2}}\sum^{[T/\Delta_n]}_{i=1}
E_{t_{i-1}} ( \Delta^n_i N_1 ) I_i
% ( |\Delta^n_i Y| \le\alpha\Delta_n^{\varpi})
\le C\Delta_n^{\varpi/2+1/4}\sum^{[T/\Delta_n]}_{i=1} E_{t_{i-1}} |
\Delta^n_i N_1 |.   \label{stable}\hspace*{-35pt}
\end{eqnarray}
%
% Taking expectation on both sides of~\ref{stable},
By Cauchy--Schwarz and Jensen's inequalities, the orthogonality of the
martingale increments, the expectation of the right-hand side in
(\ref{stable}) is
%
%e7.21 ###
\begin{eqnarray}
&\le& C\Delta_n^{\varpi/2+1/4}E\Biggl(\sum^{[T/\Delta_n]}_{i=1}
\sqrt{E_{t_{i-1}} ( \Delta^n_i N_1 )^2}\Biggr) \nonumber\\
&\le&
C\Delta_n^{\varpi/2+1/4}\frac{T}{\Delta_n}
\sqrt{\frac{\Delta_n}{T} E\Biggl(\sum^{[T/\Delta_n]}_{i=1}
E_{t_{i-1}} ( \Delta^n_i N_1 )^2\Biggr)}\\
&\le& C\Delta_n^{({\varpi-1/2)}/{2}}\sqrt{E(N_{1,T}-N_{1,0})^2}.\nonumber
\end{eqnarray}
Since $\varpi>1/2$, (\ref{stablecon}) holds. Similarly, we can
deduce that
%
%e7.22 ###
\begin{equation}
(k\Delta_n)^{{(\varpi-3/2)}/{2}}\bigl(\tilde{U}(k\Delta_n)-\sigma
_0^2\bigr)\rightarrow
\sigma_0 z_2 \label{clt2}
\end{equation}
$\mathcal{F}^{Y}$-stably. Finally in view of (\ref{clt1}) and
(\ref{clt2}), and by virtue of Lemma~\ref{le2} and~(\ref{bias}), to complete
the proof, it suffices to show that
%
%e7.23 ###
\begin{eqnarray} \label{corr}
&&\Delta_n^{3/2-\varpi}\sum^{[T/k\Delta_n]}_{i=1}E_{t_{i-1}}\Biggl(I
\bigl(|\Delta^n_{i,
k}Y|\le
\alpha(k\Delta_n)^{\varpi}\bigr)\sum^{i+k-1}_{j=i}I(|\Delta
^n_jY|\le
\alpha\Delta_n^{\varpi})\Biggr)
\nonumber\hspace*{-35pt}
\\[-8pt]
\\[-8pt]
\nonumber
&&\quad\rightarrow^P 0,\hspace*{-35pt}
\end{eqnarray}
where $\Delta^n_{i, k}Y=\sum_{j=i}^{i+k-1}\Delta^n_{j}Y$. To this
end, we give an estimate of the summands in (\ref{corr}). Let
$\Delta^{n,-j}_{i, k}Y=\Delta^n_{i, k}Y-\Delta^n_jY$, $\Delta^{n,
j-}_{i, k}=\sum^{j-1}_{l=i}\Delta_{l}^nY$ and $\Delta^{n, j+}_{i,
k}=\sum^{i+k-1}_{l=j+1}\Delta_{l}^nY$, for $i\le j\le i+k-1$. We
make the convention that $\Delta^{n, i-}_{i, k}=\Delta^{n,
(i+k-1)+}_{i, k}=0$. Then there exists a constant $C$ such that
\begin{eqnarray*}
&&\{|\Delta_{i,k}^nY|\le
\alpha(k\Delta_n)^{\varpi}\}\cap\{|\Delta_{j}^nY|\le
\alpha\Delta_n^{\varpi}\}\\
&&\qquad\subset\{|\Delta_{i,k}^{n, -j}Y|\le
C\Delta_n^{\varpi}\}\cap\{|\Delta_{j}^nY|\le
\alpha\Delta_n^{\varpi}\},
\end{eqnarray*}
%
%%
%C\Delta_n^{\varpi}\}\cup\{|\Delta_{i,k}^{n, j+}Y|\le
%C\Delta_n^{\varpi}\}\}\cap\{|\Delta_{j}^nY|\le
%%
%( \{-C\Delta_n^{\varpi}\le\Delta^{n, j+}_{i, k}Y+\Delta^{n,
%j-}_{i, k}Y\le C\Delta_n^{\varpi}\}\cap\{|\Delta^n_j Y|\le
%%
%}
and consequently, in view of $k=2$, we have
\begin{eqnarray}\label{ecorr}
& & E_{t_{i-1}}I\bigl(|\Delta^n_{i,2}Y|\le
\alpha(2\Delta_n)^{\varpi}\bigr)I(|\Delta^n_jY|\le
\alpha\Delta_n^{\varpi}) \nonumber\\
&&\qquad\le E_{t_{i-1}}[I(|\Delta^{n, j-}_{i,2}Y| \le
C\Delta_n^{\varpi})E_{t_{j-1}}I(|\Delta^n_jY|\le
\alpha\Delta_n^{\varpi})]
\nonumber
\\[-8pt]
\\[-8pt]
\nonumber
&&\qquad\quad{} + E_{t_{i-1}}[I(|\Delta^n_jY|\le
\alpha\Delta_n^{\varpi})E_{t_j}I(|\Delta^{n,
j+}_{i,2}Y|\le C\Delta_n^{\varpi})] \nonumber\\
%&& + E_{t_{i-1}}[I(|\Delta^n_jY|\le
%)] \nonumber\\
&&\qquad\le C\Delta_n^{2\varpi-1} \qquad  (\mbox{by Lemma }\ref{le2}).\nonumber
\end{eqnarray}
Substituting (\ref{ecorr}) into the left-hand side of (\ref{corr}),
we deduce that the left-hand side of (\ref{corr}) is less than
$C\Delta_n^{\varpi-1/2}$. Since $\varpi>1/2$, (\ref{corr}) is
proved.
\end{pf}

%s7.2 ###
\subsection{\texorpdfstring{Proof of Theorem \protect\ref{th2}}{Proof of Theorem 2}}

We start with the proof of the following equation which is implied
by Lemmas~\ref{lem3} and~\ref{lem4}:
%
%e7.24 ###
\begin{equation}
\Delta_n^{1+(1/\beta-\varpi)\wedge0} U(\Delta_n) \longrightarrow^P
2\alpha C_\beta\qquad  \mbox{under }  H_1, \label{prelln1}
\end{equation}
where $C_\beta$ is some constant. Then Theorem~\ref{th2} is a straight
consequence of~(\ref{prelln1}), since now $\tilde V_n \rightarrow^P
2^{1+1/\beta-\varpi}>2^{3/2-\varpi}$ and $\tilde{C}\rightarrow^P
2^{3/2-\varpi}$.

Now $X$ vanishes to a deterministic drift satisfying
$dX_t=b(X_t)\,dt$. For simplicity, we assume that
$\varepsilon^-=\varepsilon^+=:\varepsilon$. Then $Y$ admits the following
decomposition:
\begin{eqnarray*}
Y_t &=& X_t+\int^t_0\int_{|x|\le\varepsilon}x(\mu-\nu)(dx,
ds)+\int^t_0\int_{|x|>\varepsilon}x\mu(dx,
ds)\\
&&{}-\int^t_0\int_{\varepsilon<|x|\le
1}xF_s^{\prime\prime}(dx)\,ds \\
&:=& X_t+J_{1, t}+J_{2, t}+J_{3, t}.
\end{eqnarray*}

%independent L\'{e}vy process: $J_1=J_{1, 1}+J_{1, 2}$, where $J_{1,
%1}$ and $J_{1, 2}$ are associated with the L\'{e}vy densities,
%$\frac{a^{(+)}I(0<x<\varepsilon)+a^{(-)}I(-\varepsilon<x<0)}{|x|^{1+\beta}}$
%and
%$\frac{f(x)|x|^{\gamma}(a^{(+)}I(0<x<\varepsilon)+a^{(-)}I(-\varepsilon
%<x<0))}{|x|^{1+\beta}}$,
%respectively.}

The next lemma reveals that the count of small increments has almost
nothing to do with the large jumps.
\begin{lemma} \label{lem3} Under the conditions in Theorem~\ref{th2},
\[
\bigl|P_{t_{i-1}}(|\Delta^n_iY|\le
\alpha\Delta_n^{\varpi})-P_{t_{i-1}}\bigl(|\Delta^n_i(Y-J_2)|\le
\alpha\Delta_n^{\varpi}\bigr)\bigr|\le C\Delta_n.
\]
\end{lemma}
\begin{pf}
Let $M_t=\sum_{0\le s\le t}I(|\Delta_sY|> \varepsilon)$.
Then $M$ is a Poisson counting process with $\omega$ wise time
dependent intensity function
$\int_{|x|>\varepsilon}F_s^{\prime\prime}(dx)$ and
%
%e7.25 ###
\begin{equation}\qquad
P_{t_{i-1}} ( \Delta^n_i M \geq1 )\le
1-\exp\biggl(-\int_{t_{i-1}}^{t_i}\int_{|x|>\varepsilon}F_s^{\prime\prime
}(dx)\,ds\biggr)\le
C\Delta_n. \label{hasj}
\end{equation}
Notice that on $\Delta^n_i M=0$, $\Delta^n_iY=\Delta^n_i(Y-J_2)$, so
the difference within the absolute value sign is
%
%e7.26 ###
\begin{eqnarray}\label{noj}
E_{t_{i-1}}\bigl[I(|\Delta^n_iY|\le\alpha\Delta_n^{\varpi}
)-I\bigl(|\Delta^n_i(Y-J_2)|\le
\alpha\Delta_n^{\varpi}\bigr); \Delta^n_iM=0  \mbox{ or } \geq1 \bigr]
\nonumber\hspace*{-40pt}
\\[-8pt]
\\[-8pt]
\nonumber
=E_{t_{i-1}}\bigl[I(|\Delta^n_iY|\le\alpha\Delta_n^{\varpi}
)-I\bigl(|\Delta^n_i(Y-J_2)|\le
\alpha\Delta_n^{\varpi}\bigr); \Delta^n_iM \geq1 \bigr].\hspace*{-40pt}
\end{eqnarray}
Lemma~\ref{lem3} is a consequence of (\ref{hasj}) and (\ref{noj}).
\end{pf}

\begin{lemma}\label{lem4} Under Assumption~\ref{assump4},
\[
P_{t_{i-1}}\bigl(|\Delta^n_i(X+J_1+J_3)|\leq
\alpha\Delta_n^{\varpi}\bigr)=C_{\beta}\Delta_n^{(\varpi-1/\beta)\wedge
0}+o_P(\Delta_n^{\varpi-1/\beta}).
\]
\end{lemma}

\begin{pf}
Let
$\tilde{l}_i=-\alpha\Delta_n^{\varpi-1/\beta}-\Delta_n^{-1/\beta}(\Delta
^n_iX+\Delta^n_iJ_3)$
and
$\tilde{u}_i=\alpha\Delta_n^{\varpi-1/\beta}-\Delta_n^{-1/\beta}(\Delta
^n_iX+\Delta^n_iJ_3)$.
The required probability is equal to
%
%e7.27 ###
\begin{equation}
P_{t_{i-1}}(\tilde{l}_i\leq
\Delta_n^{-1/\beta}\Delta^n_iJ_1\leq\tilde{u}_i ).
\label{prob}
\end{equation}

Now we prove the lemma in two cases: (i) $\beta\geq1$ and (ii) $\beta< 1$.

\textit{Case} (i): $\beta>1$. By the L\'{e}vy--Khintchine formula,
\[
E_{t_{i-1}}\exp(i\theta\Delta_n^{-1/\beta}\Delta^n_iJ_1)=\exp(\Delta
_n\psi(\Delta_n^{-1/\beta}\theta)),
\]
where $\psi(u)=\int_R\{ \exp(iuy)-1-iuyI(|y|\leq1) \}
F^{\prime}(dy).$ By a change of variable, we have
%
%e7.28 ###
\begin{eqnarray}
\qquad\psi(\Delta_n^{-1/\beta}\theta)&=&\int_R\bigl(\exp(i\theta
z)-1-i\theta zI(|z|\leq
1)\bigr)F^{\prime}(\Delta_n^{1/\beta}z)\Delta_n^{1/\beta}\,dz
\nonumber
\\[-8pt]
\\[-8pt]
\nonumber
& &{} -\int_Ri\theta z I(1\leq|z|\leq
\Delta_n^{-1/\beta})F^{\prime}(\Delta_n^{1/\beta}z)\Delta_n^{1/\beta}\,dz.
\end{eqnarray}
Hence,
\[
\Delta_n
F^{\prime}(\Delta_n^{1/\beta}z)\Delta_n^{1/\beta}\rightarrow
\frac{1}{|z|^{1+\beta}}\bigl(a^{(+)}I(z>0)+a^{(-)}I(z<0)\bigr):=\tilde{\nu}(z).
\]
By the dominant convergence theorem,
$E\exp(i\theta\Delta_n^{-1/\beta}\Delta^n_iJ_1)$ converges to
%
%e7.29 ###
\begin{equation}
% \rightarrow
\int_R\bigl(\exp(i\theta z)-1-i\theta zI(|z|\leq
1)\bigr)\tilde{\nu}(z)\,dz+i\theta/(\beta-1)\bigl(a^{(+)}+a^{(-)}\bigr).\label{char}\hspace*{-35pt}
\end{equation}
Therefore we have $\Delta_n^{-1/\beta}\Delta^n_iJ_1$ converges in
distribution to a stable random variable. Since
$(\Delta^n_iX+\Delta^n_iJ_3)\Delta_n^{-1/\beta}=o(1)$, by
(\ref{prob}) and Assumption~\ref{assump4}, the lemma is proved in
this case.

\textit{Case} (ii): $\beta<1$. In this case, we can further decompose
$J_1$ as follows:
\[
J_1=-\int^t_0\int_{|x|\leq\varepsilon}x\nu(dx,
ds)+\int^t_0\int_{|x|\leq\varepsilon}x\mu(dx, ds):=J_{11}+J_{12}.
\]
Then the required probability in (\ref{prob}) could be rewritten as
%
%e7.30 ###
\begin{equation}
P_{t_{i-1}}(\tilde{l}_i-\Delta_n^{-1/\beta}\Delta^n_iJ_{11}\leq
\Delta_n^{-1/\beta}\Delta^n_iJ_{12}\leq
\tilde{u}_i-\Delta_n^{-1/\beta}\Delta^n_iJ_{12}).
\label{case2}\hspace*{-35pt}
\end{equation}
By similar
calculation to (\ref{char}), one gets
$\Delta_n^{-1/\beta}\Delta^n_iJ_{12}$ converges to a stable random
variable. First, consider the case where $\varpi>1/\beta$. Now by
(\ref{case2}) and Assumption~\ref{assump4}, the lemma is obtained
straightforwardly. Second, if $\varpi\leq1/\beta$, by~(\ref{case2}), the required probability is asymptotically a
constant.
\end{pf}

By Lemmas~\ref{lem3} and~\ref{lem4},
\[
\Delta_n^{1+(1/\beta-\varpi)\wedge
0}\sum^n_{i=1}E_{t_{i-1}}I(|\Delta^n_iY|\leq
\alpha\Delta_n^{\varpi})\rightarrow^P C_{\beta}2\alpha,
\]
which implies that the conditional variance goes to zero in
probability, since
\[
\sum^n_{i=1}\Delta_n^{2(1+(1/\beta-\varpi)\wedge
0)}E_{t_{i-1}}I^2(|\Delta^n_iY|\leq
\alpha\Delta_n^{\varpi})\rightarrow^P 0.
\]
Therefore, a direct use of Lenglart's inequality yields
(\ref{prelln1}).

\end{appendix}

% imsref loaded by akundreckaite, 2012-03-26 08:19:24

%suskaldyti doi

%
%

\printaddresses

\end{document}